\documentclass[11pt]{article}

\usepackage[margin=1in]{geometry}
\usepackage{amssymb}
\usepackage{amsfonts}
\usepackage{multirow}
\usepackage{amsthm}
\usepackage{graphicx}
\usepackage{enumerate}
\usepackage{dsfont}
\usepackage{amsmath,bm}
\usepackage{mathtools}
\usepackage{amsbsy}
\usepackage{authblk}
\usepackage{enumitem}
\usepackage{nameref}
\usepackage{xcolor}
\usepackage{hyperref}
\usepackage{bm}
\usepackage{natbib}
\usepackage{placeins}
\usepackage{booktabs}
\usepackage{subcaption}
\usepackage{caption}
\usepackage{threeparttable}
\usepackage{float}

\newtheorem{lemma}{Lemma}

\theoremstyle{definition}

\newtheorem{remark}{Remark}

\newtheorem{thm}{Theorem}[section]
\newtheorem{assumption}{Assumption}

\newcommand{\Cov}{{\rm Cov}}

\newcommand{\Var}{{\rm Var}}

\newcommand{\E}{\mathbb{E}}

\newcommand{\bzero}{\mathbf{0}}

\newcommand{\bh}{\mathbf{h}}

\newcommand{\bj}{\mathbf{j}}
\newcommand{\bk}{\mathbf{k}}

\newcommand{\m}{\mathbf{m}}

\newcommand{\bs}{\mathbf{s}}
\newcommand{\bt}{\mathbf{t}}

\newcommand{\bomega}{\boldsymbol{\omega}}

\def\T{^{\top}}

\newcommand{\B}{\mathcal{B}}
\newcommand{\D}{\mathcal{D}}

\newcommand{\I}{\mathcal{I}}
\newcommand{\J}{\mathcal{J}}

\newcommand{\bbr}{\mathbb{R}}

\newcommand{\te}{\theta}

\newcommand{\beq}{\begin{equation}}
\newcommand{\eeq}{\end{equation}}
\newcommand{\bea}{\begin{eqnarray}}
\newcommand{\eea}{\end{eqnarray}}
\newcommand{\beas}{\begin{eqnarray*}}
\newcommand{\eeas}{\end{eqnarray*}}

\newcommand{\keywords}[1]{%
\par\medskip
\noindent\textbf{Keywords:} #1 }

\title{\bfseries
Frequency Domain Resampling for Gridded Spatial Data
}

\author[1]{Souvick Bera\thanks{\texttt{berasouvick@mines.edu}}}
\author[2]{Daniel J. Nordman\thanks{\texttt{dnordman@iastate.edu}}}
\author[1]{Soutir Bandyopadhyay\thanks{\texttt{sbandyopadhyay@mines.edu}}}

\affil[1]{Department of Applied Mathematics and Statistics\\
Colorado School of Mines\\
Golden, Colorado 80401, USA}
\affil[2]{Department of Statistics\\ 
Iowa State University\\
Ames, IA 50011, USA}

\date{}      

\begin{document}

\maketitle

\begin{abstract}
In frequency domain analysis for spatial data, spectral averages based on the periodogram often play an important role in understanding spatial covariance structure, but also have complicated sampling distributions due to complex variances from aggregated periodograms. In order to non-parametrically approximate these sampling distributions for purposes of inference, resampling  can be useful, but previous developments in spatial bootstrap have faced challenges in the scope of their validity, specifically  due to issues in  capturing the complex variances of spatial spectral averages. As a consequence, existing  frequency domain bootstraps for spatial data are highly restricted in application to only special processes (e.g. Gaussian) or certain spatial statistics. To address this limitation and  to approximate a wide range of  spatial spectral averages, we propose 
a practical hybrid-resampling approach that combines two different resampling techniques in the forms of spatial subsampling and spatial bootstrap.  
 Subsampling helps to capture the variance of spectral averages while bootstrap captures the distributional shape. The hybrid resampling procedure can then accurately quantify uncertainty in spectral inference under mild spatial assumptions.  Moreover, compared to the more studied time series setting, this work fills a gap in the theory of subsampling/bootstrap for spatial data regarding spectral average statistics.
\end{abstract}


\section{Introduction}
\label{sec:Introduction}

In recent years, there has been a marked increase in research on the analysis of spatial data using the frequency domain approach (see, for example,~\citealp{band-etal-2017spectral, bandy2010, bandyopadhyay2015, band-rao-2016,  fuentes2006, fuentes2007, hall1994, im2007, matsuda2009, Van2017non, van2020}). As a benefit,  analysis in the frequency domain allows inference about covariance structures through a data transformation (i.e., a Fourier or periodogram-based transform), without the need for a full probability model for spatial data. Resampling methods, such as subsampling and bootstrap, for spatial data have gained popularity in recent decades because these approaches can often allow for uncertainty quantification and distributional approximations for complicated spatial statistics in an nonparametric or model-free fashion; see \cite{davison1997}, \cite{sherman1994} and \cite{sherman1998}.  In a recent study, \cite{ng2021} introduced a type of frequency domain bootstrap (FDB) for Gaussian spatial data on  grid. Although this work offers a significant contribution to nonparametric spatial approximations,  the resampling methodology comes with certain hard limitations in application. Firstly, the FDB of  \cite{ng2021}  is generally not valid for approximating many spatial spectral averages of interest, because the latter have complex  variance structures that this FDB proposal cannot correctly estimate.  In particular, the challenge is that spectral mean statistics have variances that depend both on spatial covariances as well as higher order (i.e., fourth) process cumulants, which arise due to  covariances between the spatial periodogram ordinates.  Technically, it is the variance contributions related to higher order cumulants that are generally missed or ignored in the FDB approach of \cite{ng2021}, unless the spatial data are Gaussian (i.e., so that these variance components become zero).  Consequently and secondly,   the spatial bootstrap results  of \cite{ng2021} are established under assumptions of Gaussianity, which can be stringent and may  not  suitably reflect the distributional nature of spectral  average statistics in general practice.  These aspects motivate a need to modify the spatial FDB of \cite{ng2021} to handle  a broader range of spectral statistics with application potentially \textit{non-Gaussian} spatial processes.

 To provide improved and more universally valid distributional approximations with spatial spectral statistics, we introduce a resampling technique called spatial \textit{Hybrid Frequency Domain Bootstrap} (HFDB) which aims to merge two different resampling schemes in the forms of subsampling and bootstrap into one overall approach in the frequency domain.
 The idea that spatial subsampling can be used to correctly estimate  the  variances  of complicated spectral averages, while bootstrap can be used to re-create the shape of the sampling distributions. Essentially then, subsampling plays a role in  appropriately modifying the spreads of  FDB distributional approximations, where this bootstrap would otherwise be invalid without such scaling adjustments. 
 The advantage of the proposed HFDB method  
is that this can accurately capture the uncertainty in spatial spectral statistics for calibrating tests or confidence intervals without requiring  stringent conditions on the spatial processes or assumptions of Gaussianity; in particular, the conditions needed to validate HFDB basically amount to mild assumptions for ensuring that spectral averages have limit distributions at all.  This work then intends to close some gaps in the scope of applying resampling approximations for spatial data, which is practically valuable.   Additionally, our work also bridges some gaps in the theory of subsampling and bootstrap for spatial data, and extends some recent  theoretical findings for spectral inference with time series to the spatial data setting.

We conclude this section by briefly overviewing FDB for time series and providing connections to the spatial bootstrap here. In time series analysis, FDB has  received much attention and  interest toward approximating the distributions of spectral statistics derived from the periodogram (cf.~\citealp{kreiss2012a, kreiss2012b, lahiri2003, politis2019}). The key concept behind this method is to mimic a type of asymptotic independence that exists in the periodogram  in order to generate versions of  periodogram ordinates in the bootstrap world by independent resampling (cf.~\citealp{dahlhaus1996, jentsch2010, kreiss2003, kreiss2012a, kreiss2023, kreiss2011}).
In the last three decades, there has been a progression of FDB methods, often relying on specific assumptions about the underlying time series (e.g., linear process) or differing cases of spectral mean parameters. Recently, a major advancement was achieved by \cite{meyer2020}  who introduced an innovative bootstrap method, called the hybrid periodogram bootstrap (HPB), that represents the state-of-the-art approach for time series inference concerning spectral means.  
Related to this, \cite{yu2023} recently studied subsampling  for use in conjunction with HPB, where subsampling helps to justify the latter bootstrap approach   under weaker assumptions than \cite{meyer2020} and so  HPB extends to a wide range of application with various time processes.  Unlike bootstraps that aim to re-create statistics at the same data level as the original sample, subsampling is a different resampling approach that generates smaller-scale copies of statistic; because of this distinction, subsampling often applies under weaker conditions than bootstrap, though bootstrap can provide more accurate distributional approximations when valid (cf.~\citealp{politis1999}).   The spatial FDB proposed by \cite{ng2021}
is the spatial equivalent of one of the original FBD methods for time series    
(cf.~\citealp{dahlhaus1996, jentsch2010, kreiss2003}), while the proposed spatial HFDB here intends to be a spatial analog of the HPB for time series as studied in 
\cite{meyer2020}   when combined with a subsampling extension as in \cite{yu2023}.  This development  
is non-trivial for spatial data and, similar to the time series setting, is crucial for  avoiding restrictive
process and moment conditions for dependent data.     Finally, it is worth mentioning that, similar to  \cite{ng2021}, we focus our resampling presentation on spatial observations lying on a lattice in two-dimensional space, though findings can be extended to gridded data in more general spatial sampling dimensions $d\geq 2$. 

The remainder of the paper is organized as follows. Section \ref{sec:Spectral mean statistics} introduces the distributions of spectral mean statistics, which form the basis for inference, along with the associated framework. Section~\ref{sec:sub-freq} presents a spatial subsampling framework in the frequency domain and establishes its consistency, forming the foundation for the proposed bootstrap procedure. A detailed discussion of the hybrid bootstrap method for spatial data is provided in Section~\ref{sec:HPB}. Section~\ref{sec:sim} presents numerical studies to evaluate the accuracy of the proposed procedure across several inference problems. Finally, Section \ref{sec:conclusion} concludes. Additional technical details and simulation results can be found in \nameref{sec:appendix} and \nameref{sec:supple}, respectively.

\section{Distributions of Spectral Mean Statistics}
\label{sec:Spectral mean statistics}
\subsection{Spatial process and sampling design}
\label{sec:framework}
Throughout this paper, we follow spatial sampling framework similar to \cite{bandy2010} and \cite{ng2021}.   Let $\{Z(\bs)\} : \bs \in\mathbb{Z}^2\}$ denote a real-valued second-order stationary process located on a regular integer grid  $ \mathbb{Z}^2$.  In this fashion, potential observations lie on a spatial lattice with a constant separation  in each coordinate direction (which we take to be 1 though other scaling is possible as well).  We assume that the random process $Z(\bs)$ is observed at $n\equiv n_1n_2$  sampling sites 
defined by 
\begin{eqnarray*}
    \{\bs_1,\ldots,\bs_{n}\} \equiv \{\bs \in \mathbb{Z}^2 : \bs \in \D_n\} = \D_n\cap \mathbb{Z}^2
\end{eqnarray*}
given by those locations on the grid $\mathbb{Z}^2$ that lie within a rectangular   sampling region   $\mathcal{D}_n\equiv [1, n_1]\times [1, n_2]$  for integers $n_1,n_2 \geq 1$.  
 For spatial data analysis, the above sampling framework  corresponds to a so-called   pure increasing domain for developing spatial results (cf.~\citealp{cressie1993}). That is, in this  scheme, more spatial observations are available (i.e., $n\to \infty$) as spatial sampling region expands 
 in size in each direction (i.e., $n_i \rightarrow \infty$ for $i=1,2$).

\subsection{Spectral mean parameters}
\label{sec:spectral mean parameter}

Suppose that the second order stationary  spatial process $\{Z(\bs):\bs \in \mathbb{Z}^2\}$ has  a spectral density denoted as $f(\bomega) : \Pi ^2 \longrightarrow [0, \infty)$ where $f(\bomega)\equiv  (2\pi)^{-2}
  \sum_{\bh \in \mathbb{Z}^2}    \Cov(Z(\bzero), Z(\bh))  \exp(-\iota  \bh^{\top} \bomega)$, for $\iota\equiv \sqrt{-1}$ and for  
$\bomega \in \Pi^2  \equiv [-\pi,\pi]^2$.
Then, the target spatial parameter for inference is a \textit{spectral mean} parameter defined as an integral
\begin{equation}
\label{eq:spectral-mean}
M(\psi) \equiv \int_{\Pi^2}\psi(\bomega)f(\bomega)d \bomega,
\end{equation}
involving $f$ and some chosen/given function $\psi:\Pi^2\longrightarrow \mathbb{R}$ of bounded variation. Spectral mean parameters are common in the spatial data analysis for studying the spatial covariance structure of the data with an unknown distribution.  Such quantities have a well-established history for both time series and spatial data (cf.~\citealp{politis2019}, for  details), dating as far back as \cite{parzen1957}, and many spatial parameters of interest can be expressed in the form of a spectral mean. Standard examples of spectral means include auto-covariances
and spectral distributions, as described below, depending on the choice of $\psi$ above.

Next we give some examples of spectral mean parameters arising in frequency domain analysis (cf.~\citealp{bandyopadhyay2015, dahlhaus1985, dahlhaus1996, parzen1957, subba2018,Van2017non,van2020} and the references therein).
\subsubsection*{Examples of spectral mean parameters}
\label{sec:examples}

\noindent\textit{1. Covariance function.} For a lag vector $\bh\in\mathbb{Z}^2$, if we consider the function $\psi(\bomega) =\cos(\bh^{\top}\bomega)$ in \ref{eq:spectral-mean}, this leads to the autocovariance function $\gamma(\bh) \equiv \Cov(Z(\bzero), Z(\bh))$, which can be viewed as a spectral mean $M(\psi)$. 
Note that, for testing the special hypothesis $H_{0}: Z(\bs)$ is white noise (constant $f$), a test based on the process auto-covariances can be applied, similar to some Portmanteau tests (cf.~\citealp{li1986,ljung1978}) using the fact that the spectral means $M(\psi)$
are zero under this $H_{0}$ for the function choice $\psi(\bomega) = [\cos(\bh_{1}{\T}\bomega),\ldots, (\cos(\bh_{k}{\T}\bomega)]{\T}$ for some integer lags $\bh_1,\ldots,\bh_k$.

\vskip0.5em

\noindent\textit{2. Spectral distribution function.} For a vector $\bt\in\mathbb{R}^2$, the spectral distribution function is defined as $F(\bt) \equiv \int_{\mathbb{R}^2}\mathbb{I}_{(-\infty, \bt]}(\bomega)f(\bomega) d\bomega$, which corresponds to a spectral mean parameter defined by 
$\psi(\bomega) = \mathbb{I}_{(-\infty, \bt]}(\bomega)$ in (\ref{eq:spectral-mean}) where $\mathbb{I}(\cdot)$ is the indicator function and $(-\infty, \bt] \equiv (-\infty, t_{1}] \times (-\infty, t_{2}]$. The spectral distribution function plays a  role in determining the smoothness of the sample paths of the spatial process $Z(\cdot)$ (cf.~\citealp{stein1999}). 

\vskip0.5em

\noindent\textit{3. Assessment of spatial covariance structures.} When analyzing spatial data, it can be useful to assess the spatial covariance structures in a nonparametric manner, without requiring specific distributional assumptions about the data. General testing methods can be developed  of evaluating different hypotheses about spatial covariances, such as tests of isotropy or separability, by examining $M(\psi)$ with appropriate choices of $\psi(\bomega)$ functions. For more details,   refer to \cite{Van2017non, van2020}.

\vskip0.5em

\noindent\textit{4. Whittle Estimation.} Suppose $\left\{f_\theta\right\}$ represents a parametric family of spectral densities (indexed by $\theta$). Whittle estimation aims to identify the closest member $f_\theta$ to the true density $f$ (cf.~\citealp{taniguchi1979}). Assuming real-valued $\theta$ for illustration, the solution to a spectral mean   $M(\psi)=$ 0 identifies the appropriate $f_\theta$ under certain conditions, using a function $\psi(\bomega)  =[1-$ $\left.f_\theta(\bomega) / f(\bomega)\right] \frac{d}{d \theta} f_\theta^{-1}(\bomega), \bomega \in \mathbb{R}^2$, defined by the derivative of $f_\theta^{-1}(\bomega) \equiv 1 / f_\theta(\bomega)$.

\vskip0.5em

\noindent\textit{5. Goodness-of-fit tests.} There has been increasing interest in frequency domain based tests to assess model adequacy (cf.~\citealp{crujeiras2010a, van2020, weller2020}). Consider a test
involving a simple null hypothesis $H_{0} : f=f_{0}$ against an alternative $H_{0} : f\neq f_{0}$ for some candidate spectral density $f_0$. One immediate test for $H_0$ is based on the function $\psi(\bomega)=1/f_{0}(\bomega)$ with $M(\psi)$ a constant   under this simple $H_0$.
To test the composite hypothesis $H_{0}: f\in \mathcal{F}$ for a specified parametric class $\mathcal{F}$, several frequency domain tests have been proposed in time series (cf.~\citealp{beran1992}, \citealp{milhoj1981},  \citealp{paparoditis2000}, \citealp{nordman2006}). These tests can use Whittle estimation to choose the ``best'' fitting model from $\mathcal{F}$ and then compare the fitted density to the periodogram across all ordinates. One can adapt   strategies for goodness-of-fit tests with time series (cf.~\citealp{nordman2006}), which combine 
aspects of model fitting and model comparison into a choice of  estimating function $\psi$, to spatial processes. 

\vskip0.5em

\noindent\textit{6. Variogram model fitting.} A popular approach  to fitting a parametric variogram model to  spatial data is through the method of least squares; see \cite{cressie1993}. Let $\{2\nu(\cdot;\theta): \theta \in \Theta\}$, $\Theta\subset \bbr^p$ denote a class of valid variogram models for the true variogram $ 2\nu(\bh)  \equiv \Var (Z(\bh) - Z(\bzero))$, $\bh\in \mathbb{Z}^2$ of the spatial process. Let $2\widehat{\nu}_n(\bh)$ denote the sample variogram at lag $\bh$ based on $Z(\bs_1),\ldots,Z(\bs_n)$ (cf.~Chapter 2, \citealp{cressie1993}). Then one can  fit the variogram model by  estimating   the parameter $\te$ that  minimizes the criterion: $\widehat{\theta}_{n}\equiv  \mbox{argmin} \Big \{ \sum_{i=1}^m  \Big( 2\widehat{\nu}_n(\bh_i)
- 2 \nu(\bh_i;\te)\Big)^2 : \te \in \Theta\Big\}$
for a given set of lags $\bh_1,\ldots, \bh_m$. Expressing the variogram in terms of the spectral density function, we get an equivalent  spectral estimating  equation for identifying $\theta$: $M(\psi)=0$ for 
$$
M(\psi)\equiv \int_{\bbr^2} \Big[\sum_{i=1}^m  \Big\{1 - \cos (\bh_i{\T}\bomega)-\nu(\bh;\theta)\Big\}
                \nabla [2\nu(\bh_i;\theta)]\Big]f(\bomega) d\bomega,
$$
using $\psi_\theta(\bomega) =\sum_{i=1}^m  \Big\{1 - \cos (\bh_i{\T}\bomega)-\nu(\bh;\theta)\Big\}
                \nabla [2\nu(\bh_i;\theta)]$.

\subsection{Spectral mean statistics and sampling distributions}
\label{sec:spec-distn}

From the spatial process $Z(\cdot)$ observed on the sampling region $\mathcal{D}_n$, where $n = n_{1}n_{2}$ is the number of observations,   we define $\I_{n}(\bomega)$ to be the periodogram at a frequency $\bomega \in \Pi^{2}$ as
$$ \I_{n}(\bomega)\equiv  (2\pi)^{-2} n^{-1}\left|\sum_{s_1=1}^{n_1}\sum_{s_2=1}^{n_2} Z(  \bs) \exp{(-\iota{ \bs^{\T} \bomega )}}\right|^2,\ \bs \equiv  (s_1, s_2)^{\T},\ \mbox{where}\ \iota\equiv \sqrt{-1}.$$    
Using $\I_{n}(\cdot)$ in place of $f(\cdot)$ in  (\ref{eq:spectral-mean}), a standard estimator of the spectral mean parameter $M(\psi)$ is then the \textit{spectral mean statistic}, or {\it spectral average}, given in Riemann sum form by
\begin{equation*}
\widehat{M}_n(\psi) \equiv (2\pi)^2 n^{-1} \sum_{\bj\in \J_{n}} \psi(\bomega_{\bj,n}) \I_n(\bomega_{\bj,n}),
\end{equation*}
using discrete frequencies $\bomega_{\bj,n}\equiv  (2\pi j_1/  n_1,2 \pi j_2/  n_2)$ defined by  $\bj \equiv (j_1,j_2) \in \J_{n}$ with index set $\J_{n} \equiv  \{\lfloor{-(n_1 - 1)/2}\rfloor,\ldots, n_1- \lfloor{n_1/2}\rfloor \} \times \{\lfloor{-(n_2 - 1)/2}\rfloor,\ldots, n_2- \lfloor{n_2/2}\rfloor \} \setminus \{(0,0)\}$ denoting the (nonzero) discrete Fourier frequency grid. For calibrating tests and confidence intervals, a spectral mean statistic has a large-sample normal distribution under mild conditions (cf.~\citealp{brillinger2001, dahlhaus1985}) given by

\begin{equation}
\label{eq:tphi}
    H_n(\psi)  \equiv  n^{1/2}\{\widehat{M}_{n}(\psi)-M(\psi)\} \overset{D}{\longrightarrow} \mathcal{N}(0, \sigma^2 = \sigma^2_1 + \sigma^2_2) \hspace{0.3em} \textnormal{as} \hspace{0.4em} n \rightarrow \infty,
\end{equation}
with variance components
\begin{eqnarray*} 
    \sigma^2_1 \equiv \sigma^2_1(\psi) &=&  (2\pi)^2 \int_{\Pi^{2} } \psi(\bomega)[\psi(\bomega)+\psi(-\bomega)]f^{2}(\bomega)\,d\bomega, \\
    \sigma^2_2 \equiv \sigma^2_2(\psi) &=&  (2\pi)^2 \int_{\Pi^{2} }\int_{\Pi^{2} } \psi(\bomega_{1})\psi(\bomega_{2})f_{4}(\bomega_{1},\bomega_{2},-\bomega_{2})\,d\bomega_{1}\,d\bomega_{2}
\end{eqnarray*}
where, for $\bomega_{1}, \bomega_{2}, \bomega_{3} \in \Pi^{2}$, 
$$
f_4(\bomega_{1}, \bomega_{2}, \bomega_{3}) \equiv  (2\pi)^{-3}
  \sum_{\bh_{1}, \bh_{2}, \bh_{3} \in \mathbb{Z}^2} \mbox{cum}(Z(\bzero), Z(\bh_1), Z(\bh_2), Z(\bh_3)) \exp\left(-\iota \sum_{\ell=1}^{3} \bh_\ell^{\T} \bomega_{\ell}\right),
$$
denotes the 4th-order cumulant spectral density. Given the complicated  structure in the limit distribution  from (\ref{eq:tphi}), one attractive approach is to consider   resampling for approximating the distribution of $H_n(\psi)$ nonparametrically. However, estimation  of the variance component $\sigma_2^2$, in particular, poses significant challenges for bootstrap approximations to the distribution of $H_n(\psi)$, as needed for the inference of spatial mean parameters. Specially, the spatial FDB introduced in \cite{ng2021} cannot capture this component  in general, unless the underlying process is assumed to be Gaussian in which case this variance component $\sigma_2^2 = 0$ becomes zero.  As a consequence, in order to improve bootstrap inference for spatial data,    one must first correctly estimate all variance components in  (\ref{eq:tphi}) in the resampling mechanism used.  For use in conjunction with spatial bootstrap, we next consider 
a spatial subsampling for the purpose of such variance estimation, as described in the following section. 

\section{Spatial Subsampling in Frequency Domain}
\label{sec:sub-freq}
\subsection{Spatial subsampling  variance estimators} \label{sec:SFDS}
Based on the observed spatial data $\{Z(\bs): \bs \in\mathcal{D}_n\cap  \mathbb{Z}^2\}$ lying within the  spatial sampling region $\mathcal{D}_n\equiv [1,n_1]\times[1,n_2]$, subsampling aims to create several ``smaller scale copies" of the original spatial data by using  data 
blocks  of smaller size $b_{n} \equiv n_{1}^{(b)}  n_{2}^{(b)}$  with $n_{k}^{(b)}<n_{k},\ k=1,2$ (cf.~\citealp{lahiri2003, sherman1994,sherman1996,sherman1998})  compared to  $n=n_1   n_2$,  where   the  block size $n_k^{(b)}$ and the original sample size $n_k$  in each direction   satisfy the relation $n_k^{(b)}\rightarrow \infty$ and $n_k^{- 1}n_k^{(b)}\rightarrow 0$ as $n \rightarrow \infty$.  In particular, we can define data blocks by all integer translates of a region $\mathcal{D}_b \equiv [1,n_1^{(b)}]\times [1,n_2^{(b)}]$
that lie inside the original sampling region $\mathcal{D}_n \equiv [1,n_1]\times[1,n_2]$; that is, each data block is of the form $\bj + \mathcal{D}_b$ 
for an integer $\bj\equiv (j_1,j_2)$ with components $j_k=0,\ldots, n_k-n_k^{(b)}$, $k=1,2$.  
In total, there are, say $L\equiv  (n_1-n_1^{(b)}+1)(n_2-n_2^{(b)}+1)$, such data blocks and, for simplicity, we  denote the $\ell$-data block as 
$\mathcal{B}_\ell$ for $\ell=1,\ldots, L$.   
We  then define a subsample periodogram for the $\ell$-th block as 
\begin{equation*}
    \I^{(\ell)}_{sub}(\bomega) \equiv (2\pi)^{-2}b_{n}^{-1} \hspace{.1cm} \left|\sum_{\bs \in \B_{\ell}\cap \mathbb{Z}^2} Z(\bs) \hspace{.1cm} e^ {-\iota  \bs^{\T} \bomega}\right|^2, \hspace{.2cm} \bomega \in \Pi^{2}, \hspace{.1cm} \ell = 1,\ldots,L.
\end{equation*}
The corresponding subsample analog of the mean statistic $\widehat{M}_{n}(\psi)$ can then be defined as 
\begin{equation*}
\widehat{M}^{(\ell)}_{sub}(\psi) \equiv (2\pi)^2 b_{n}^{-1} \sum_{\bj \in \J_b} \psi(\bomega_{\bj,b}) \I_{sub}^{(\ell)}(\bomega_{\bj,b}),
\end{equation*}
where  $\bomega_{\bj,b}$ and $\J_b$ are defined similar to   $\bomega_{\bj,n}$ and $\J_n$  earlier, but  based on the subsample sizes $b_n$ rather than the original data size $n$.  
These lead to a collection of subsample versions $H^{(\ell)}_{sub}(\psi)$ of the target spectral quantity $H_n(\psi) \equiv n^{1/2} \{ \widehat{M}_n(\psi) - M(\psi) \}$ 
as
\begin{equation*}
     H^{(\ell)}_{sub}(\psi) \equiv b_{n}^{1/2} \{ \widehat{M}^{(\ell)}_{sub}(\psi) - \widetilde{M}_n(\psi) \},\ \ell= 1,\ldots,L,
\end{equation*}
with the subsample mean $\widetilde{M}_n(\psi) \equiv L^{-1} \sum_{\ell=1}^{L} \widehat{M}^{(\ell)}_{sub}(\psi)$ as proxy for  the  parameter $M(\psi)$. 

Recall that spatial variance estimation is especially important here because the  variance $\sigma^2$ of spectral mean statistics from (\ref{eq:tphi}) can be quite complex or difficult to determine accurately, which is our motivation
for consideration of subsampling.  Consequently, we can define the subsampling variance estimator of $\sigma^2$ as 
\begin{equation}
\label{eqn:sigma}
    \widehat{\sigma}^2_n = \widehat{\sigma}^2_n(\psi) \equiv \int_{\mathbb{R}} z^2 d \widehat{F}_{H_n(\psi)}(z) = L^{-1} \sum_{\ell=1}^{L} b_n\{\widehat{M}^{(\ell)}_{sub}(\psi) - \widetilde{M}_n(\psi) \}^2,
\end{equation}
which corresponds to 
the sample variance of the subsample statistics $  \{ H^{(\ell)}_{sub}(\psi) \equiv b_{n}^{1/2} \{ \widehat{M}^{(\ell)}_{sub}(\psi) - \widetilde{M}_n(\psi) \}_{\ell=1}^{L}$.
Further, using subsampling, we also introduce an estimator for just the first variance component $\sigma^2_1 $  in  (\ref{eq:tphi})  as
\begin{eqnarray}
\label{eq:estvar_first}
   \widehat{\sigma}^2_{1,n} \equiv b_n^{-1} (4\pi^2)^2 \sum_{\bj\in \J_{b}} \psi(\bomega_{\bj,b})(\psi(\bomega_{\bj,b})+\psi(\bomega_{-\bj,b}))\frac{1}{L} \sum_{\ell=1}^{L}(\I^{(\ell)}_{sub}(\bomega_{\bj,b})-\widetilde{\I}_n(\bomega_{\bj,b}))^2,
\end{eqnarray}
where $\widetilde{\I}_n(\bomega_{\bj,b}) \equiv L^{-1}\sum_{\ell=1}^{L} \I^{(\ell)}_{sub}(\bomega_{\bj,b})$.
Here in (\ref{eq:estvar_first}), we are isolating a component from the overall subsampling variance estimator $\widehat{\sigma}^2_n$ in (\ref{eqn:sigma})  that arises due to    marginal variance  in 
subsample periodogram copies $\bigl\{\I^{(\ell)}_{sub}(\bomega_{\bj,b})\bigr\}_{\ell=1}^L$.    This leads to a second subsampling estimator
\begin{eqnarray}
\label{eq:estvar_second}
    \widehat{\sigma}^2_{2,n} \equiv \widehat{\sigma}^2_n-\widehat{\sigma}^2_{1,n}
\end{eqnarray}
to approximate  the second or remaining component $\sigma^2_2 $ in    (\ref{eq:tphi}).  Spatial subsampling then enables a valid estimation of both variance components $\sigma_1^2$ and $\sigma_2^2$ which contribute to the target distribution of $H_n(\psi)$ having a limit variance of $\sigma^2 \equiv \sigma_1^2 +\sigma_2^2$.

\subsection{Consistency   of spatial subsampling estimators}
\label{sec:cons_subest}
To establish the formal consistency properties for subsampling variance estimators with  spectral mean statistics, we use the following mild assumptions.  Assumptions \ref{assumption:assump1}-\ref{assumption:assump2}   are standard for ensuring that
 the target quantity $H_n(\psi)$ in (\ref{eq:tphi}) has a limit law with    finite variance $\sigma^2 $. 
\begin{assumption}
\label{assumption:assump1}
    $\{Z( \bs):\bs \in \mathbb{Z}^2\}$ is 4-th order stationary with   summable covariances
    and 4-th order cumulants in  the sense that
      $\sum_{\bs\in\mathbb{Z}^2} (1+ \|\bs\|)|\mathrm{Cov}(Z(\bm{0}),Z(\bs))| <\infty$ 
      and $ \sum_{\bs_1,\bs_2,\bs_3\in\mathbb{Z}^2}  |\mathrm{cum}(Z(\bm{0}),Z(\bs_1),z(\bs_2),z(\bs_3))| <\infty$,
    where $\|\cdot\|$ denotes the Euclidean norm.
\end{assumption}
\begin{assumption}
\label{assumption:assump2}
    The normal limit     exists in (\ref{eq:tphi}) for any $\psi$ of bounded variation.
\end{assumption}

We require  two 
additional 
Assumptions~\ref{assumption:assump3}-\ref{assumption:assump4}    below, which are related 
 to Assumptions \ref{assumption:assump1}-\ref{assumption:assump2}, but typically weaker (i.e., spatial processes fulfilling \ref{assumption:assump1}-\ref{assumption:assump2}  typically then satisfy
\ref{assumption:assump3}-\ref{assumption:assump4}; cf.~Remark \ref{rem:remark1}).
Assumptions~\ref{assumption:assump3}-\ref{assumption:assump4} entail that  dependence should decrease between certain statistics  computed on distant observational regions, which  is necessary  for the validity of subsampling estimation
(cf.~Theorem~\ref{thm:consistency-variance}) that relies on  data blocks to provide spatial replication.

For two sequences $\bj^{(n)} \equiv (j^{(n)}_1, j^{(n)}_2),  \bk^{(n)} \equiv (k^{(n)}_1, k^{(n)}_2)\in \mathbb{Z}^2$
of integer pairs,  
define  $ H_{\bj^{(n)}, n}(\psi)$ and $ H_{\bk^{(n)}, n}(\psi)$ as versions of $H_n(\psi)$  when computed from   translated spatial regions $\bj^{(n)} + \D_n$ and  $\bk^{(n)} + \D_n$, instead of the region $ \D_n \equiv [1,n_1]\times [1,n_2]$ used for  $H_n(\psi)$  in (\ref{eq:tphi}).

\begin{assumption}
\label{assumption:assump3}
For any two  $\bj^{(n)}$ and  $\bk^{(n)}$ integer sequences where  $\max_{i=1,2} (|j^{(n)}_i-k^{(n)}_i|/n_i )\rightarrow \infty$    as $n\equiv n_1  n_2 \rightarrow \infty$,    the following quantities are  asymptotically normal and independent:  
\begin{eqnarray*}
\begin{pmatrix}
  H_{\bj^{(n)}, n}(\psi) \\
  H_{\bk^{(n)}, n}(\psi)
\end{pmatrix} 
 \overset{D}{\longrightarrow}  \mathcal{N} \left(\begin{pmatrix}
     0 \\ 0
 \end{pmatrix}, \sigma^2 \left(\begin{matrix}
     1&0 \\ 0&1
 \end{matrix}\right)\right) \quad \mbox{as $n\to \infty$.}
\end{eqnarray*}

\end{assumption}
 For two integer sequences  $\bj^{(n)},\bk^{(n)} \in \mathbb{Z}^2$,
 we use one further assumption stated below that relates to the periodogram, say $\I_{\bj^{(n)}, n}(\bomega)$ and $\I_{\bk^{(n)}, n}(\bomega)$, when computed on translated sampling regions $\bj^{(n)} + \D_n$ and $\bk^{(n)} + \D_n$, respectively, for $\bomega   \in \Pi^2 \equiv [-\pi,\pi]^2$.

\begin{assumption}
\label{assumption:assump4}
 Let $\m_n  \equiv (m_{n_1}, m_{n_2}) \in\mathbb{Z}^2$ with $0 \leq  |m_{n_i}| \le n_i - \lfloor n_i/2 \rfloor, i=1,2 $ be an integer sequence such that the discrete Fourier frequencies $\bomega_{\m_n,n}$ converge to some limit $\bomega \equiv (\omega_1, \omega_2)$ with $0<|\omega_i|<\pi$, $i=1,2$, as $n \equiv n_1   n_2 \to \infty$. 
For any two  $\bj^{(n)}$ and  $\bk^{(n)}$ integer sequences where  $\max_{i=1,2} (|j^{(n)}_i-k^{(n)}_i|/n_i )\rightarrow \infty$  as $n  \rightarrow \infty$,   the  periodogram ordinates are  asymptotically exponential and independent: 
\begin{eqnarray*}
\begin{pmatrix}
    \I_{\bj^{(n)}, n}(\bomega_{\m_n,n})  \\
    \I_{\bk^{(n)}, n}(\bomega_{\m_n,n})
\end{pmatrix} 
 \overset{D}{\longrightarrow}  f(\bomega) \begin{pmatrix}
     U_1 \\ U_2
 \end{pmatrix} &  \quad \mbox{as $n\to \infty$}
\end{eqnarray*}  
\end{assumption}
 \noindent where $U_1, U_2$ denote i.i.d standard exponential random variables.

For context,
Assumption \ref{assumption:assump3} intends a general, but mild, characterization of weak spatial dependence  connected to Assumptions \ref{assumption:assump1}-\ref{assumption:assump2}. This condition states that any
copies $H_{\bj^{(n)},n}(\psi)$ and $H_{\bk^{(n)},n}(\psi)$  of the statistical quantity of interest  $H_{n}(\psi)$,  similarly defined on   size $n\equiv n_1 n_2$ regions but of a diverging distance   as $n\rightarrow \infty$,
 should be normal like $H_n(\psi)$ but also asymptotically independent; in particular,  the spatial regions $\bj^{(n)} + \D_n$ and $\bk^{(n)} + \D_n$ for defining 
$H_{\bj^{(n)},n}(\psi)$ and $H_{\bk^{(n)},n}(\psi)$  have a distance 
\begin{equation}
\label{eqn:dis}
\mathrm{dist}_n \equiv \min\Big\{ \|\bs -\bt\| : \bs \in \bj^{(n)} + \D_n,\, \bt \in \bk^{(n)} + \D_n\Big\} \geq \max_{i=1,2}\left(|j^{(n)}_i - k^{(n)}_i | - n_i
\right)  \rightarrow \infty  
\end{equation}
that grows as $n\to \infty$, because  
$i$th component of the difference $\bj^{(n)} - \bk^{(n)}$ diverges faster than the sample size $n_i$ in some direction, $i=1,2$, in Assumption~\ref{assumption:assump3}. 
For such regions,  a natural characterization of weak dependence is that   quantities $H_{\bj^{(n)},n}(\psi)$ and $H_{\bk^{(n)},n}(\psi)$ should be independent as $n\rightarrow \infty$. 
 Assumption \ref{assumption:assump4} is analogous to Assumption \ref{assumption:assump3}  in spirit,  but  weaker, saying  that spatial periodograms computed from distant regions should be asymptotically independent and  exponentially distributed up to scaling by the spectral density $f$. 
Assumptions~\ref{assumption:assump1}-\ref{assumption:assump4} 
can hold for many spatial processes, and we provide some examples next.    

As perhaps the most straightforward example, an   $m$-dependent stationary spatial process with 
$\E [Z(\bs)]^4<\infty$  
satisfies Assumptions~\ref{assumption:assump1}-\ref{assumption:assump4} (i.e., where two observations $Z(\bs)$ and $Z(\bt)$ are independent if $\|\bs-\bt\|>m$ for some distance $m>0$); in particular, the independence of statistics in Assumption~\ref{assumption:assump3}  holds by (\ref{eqn:dis}) in this case.  Examples of such processes can include Gaussian nearest-neighbor models (cf.~\citealp{datta2016nearest}).    
Another standard spatial process  where Assumptions~\ref{assumption:assump1}-\ref{assumption:assump4} hold is given by linear spatial processes $Z(\bs) = \mu+\sum_{\bk \in \mathbb{Z}^2} a_{\bk} \varepsilon(\bs -\bk)$ defined by i.i.d. (perhaps non-Gaussian) innovations $\{\varepsilon(\bk):\bk\in\mathbb{Z}^2\}$ with  $\E[\varepsilon(\bk)]=0$ and $\E[\varepsilon(\bk)]^4<\infty$ and real-valued coefficients that are summable $\sum_{\bk \in\mathbb{Z}^2} \|\bk\| |a_{\bk}|<\infty$.  These processes have a long history of consideration in statistics and econometrics (cf.~\citealp{tjostheim1978, anselin1988spatial, guyon1995random}), and   Assumptions~\ref{assumption:assump1}-\ref{assumption:assump4} can be verified by approximating such processes with tail-truncated versions $\mu+\sum_{\|\bk\|<m/2} a_{\bk} \varepsilon(\bs -\bk)$ (i.e., for suitably large $m$) that are $m$-dependent.  A further general class  of spatial processes that   fulfill Assumptions~\ref{assumption:assump1}-\ref{assumption:assump4}  include so called $\alpha$-mixing  processes, which are commonly considered in resampling developments   (cf.~\citealp{lahiri2003}); see Remark \ref{rem:remark1} for   technical details.    Non-Gaussian examples of such processes can include spatial Markov random fields, having conditionally specified distributions that may be
  beta, binary, Poisson, or more general exponential families (cg.~\citealp{besag1974,caragea2009autologistic, hardouin2008multi, kaiser2002spatial}). In particular, when such conditional model specifications satisfy  Dobrushin’s uniqueness condition (cf.~\citealp{guyon1995random}), then the Markov random field is $\alpha$-mixing  (actually, stronger $\phi$-mixing).

We can next state a formal result on the consistency of subsampling variance estimators.   

\begin{thm}
\label{thm:consistency-variance}
Suppose Assumptions \ref{assumption:assump1}-\ref{assumption:assump4} hold and the subsample 
size $b_n \equiv n_1^{(b)}   n_2^{(b)}$ satisfies  $1/n_k^{(b)} + n_k^{-1} n_k^{(b)} \rightarrow 0$ for $k=1,2$ as $n\equiv n_1  n_2\rightarrow \infty$.  Then, the  spatial subsampling estimators of variance components  are consistent:
\begin{equation*}
\widehat{\sigma}^2_{1,n} \overset{P}{\longrightarrow} \sigma^2_1,\quad 
\widehat{\sigma}^2_{2,n} \overset{P}{\longrightarrow} \sigma^2_2,\ \ \mbox{and}\ \ \ \widehat{\sigma}^2_{n} \overset{P}{\longrightarrow}\sigma^2\equiv \sigma_1^2 +\sigma_2^2 \quad \mbox{as $n\to \infty$.}
\end{equation*}
\end{thm}
The theorem above provides a device  to modify and improve spatial bootstraps 
for approximating spectral mean statistics, particularly when such bootstraps can fail to appropriately reflect the spreads of such statistics.  This then leads to the spatial hybrid bootstrap procedure proposed in the next section.

\begin{remark}
\label{rem:remark1}

For completeness, we state $\alpha$-mixing conditions for verifying  Assumptions~\ref{assumption:assump1}-\ref{assumption:assump4}. For a subset $T \subset \mathbb{Z}^2$ of spatial locations,  let $|T|$ denote the cardinality and let  $\mathcal{F}_Z(T) \equiv \sigma \langle Z(\bs): \bs \in T \rangle$  denote the $\sigma$-algebra   generated by observations $Z(\bs)$, $\bs \in T$. Define a dependence measure
  $\tilde{\alpha}(T_1,T_2) \equiv \sup \{ |P(A_1 \cap A_2) - P(A_1)P(A_2)|: A_1 \in \mathcal{F}_Z(T_1),A_2 \in \mathcal{F}_Z(T_2) \}$ and distance 
$d(T_1,T_2) \equiv \min\{\|\bs_1 -\bs_2\|: \bs_1 \in T, \bs_2\in T_2\}$ 
between two subsets $T_1,T_2 \subset \mathbb{Z}^2$.  Denote the mixing coefficient 
of the spatial process
$\alpha(a,b)  \equiv \sup \{\tilde{\alpha}(T_1,T_2) : d(T_1,T_2) \geq a, \max\{|T_1|, |T_2|\}\leq b \}$ as a function of the distance $a>0$ and size $b>0$ of two sets  of spatial locations. Assuming $\alpha(a,b) \leq \alpha_1(a) g(b)$ holds for a non-increasing function $\alpha_1(\cdot)$ and a non-decreasing function $g(\cdot)$ as in~\cite{lahiri2003central}, results in the latter (Prop 4.2, Theorem 4.2) can be used to show that  
Assumption~\ref{assumption:assump1} holds if 
$\sup_{\bs \in \mathbb{Z}^2}\E[Z(\bs)]^{4 + \delta} <\infty$ along with 
$ \alpha_1(a)=O(a^{-\tau})$ as $a\to \infty$ for  some $\delta>0$ and  $\tau> 6(4+\delta)/\delta$, while  Assumption~\ref{assumption:assump2} also holds if 
$g(b) = o( b^{(\tau-2)/8})$ as $b\to \infty$
and the directional sizes     $n_1,n_2$ of a sampling region   
are proportional (or $\max\{n_1,n_2\}/\min\{n_1,n_2\}$ is bounded). Such conditions  further imply Assumptions~\ref{assumption:assump3}-\ref{assumption:assump4}, as the   dependence between the two sets/statistics in \ref{assumption:assump3}-\ref{assumption:assump4} is measured by  $\alpha( \mathrm{dist}_n,n) = O(n_1^{-\tau} n^{(\tau-2)/8})\rightarrow 0$ as $n=n_1   n_2 \to \infty$ for   $\mathrm{dist}_n$ in (\ref{eqn:dis}).
\end{remark}

\section{Spatial Hybrid Frequency Domain Bootstrap (HFDB)}
\label{sec:HPB}
We may now describe our main resampling approach for approximating the distribution of spatial spectral mean statistics. Recall that, as explained in the Introduction, \cite{ng2021}  proposed a bootstrap for spatial data, which is generally invalid for inference about spectral averages  and requires modification for capturing the complicated variance of spectral averages.  For reference, we shall term their bootstrap method as technically being a    Frequency Domain Wild Bootstrap (FDWB).  Recently,  for the case of time series data, \cite{meyer2020} and \cite{yu2023} described a hybrid periodogram bootstrap  (HPB)
 as the most advanced  resampling scheme for spectral mean inference with time series in terms of broad validity and performance.  Our hybrid resampling approach for spatial data, termed the Hybrid Frequency Domain Bootstrap (HFDB), intends 
 types of extensions of both FDWB  and HPB methods.  That is, the spatial HFDB aims to   combine spatial bootstrap (i.e., FDWB) with spatial subsampling (Sec.~\ref{sec:SFDS})
in order to correct scaling issues in bootstrap for spatial data.  This also produces a spatial analog version of HPB from time series.     

For setting distributional approximations for spatial spectral mean statistics,  both HFDB and FDWB approaches start with a common step of mimicking the distribution of spectral statistics by a scheme of   independently resampling or re-creating  periodogram values based on an estimator $\widehat{f}_n$ of the spectral density $f$.    To this end, let $\widehat{f}_n$ represent a consistent estimator of the spectral density $f$ across spatial discrete Fourier frequencies, i.e.
\begin{equation}
\label{eq:specden-consistency}
    \max\limits_{\bj \in \J_n} \lvert \widehat{f}_n(\bomega_{\bj,n}) - f(\bomega_{\bj,n})  \rvert = o_P(1) \hspace{0.2cm} \mbox{as}\ \ n\rightarrow \infty.
\end{equation} 
 Note that all FDB approaches involve spectral density estimation, which is true for time series  (cf.~\citealp{dahlhaus1996, jentsch2010, kreiss2003, meyer2020}) and for FDWB with spatial data (\citealp{ng2021}).
We then define an initial bootstrap re-creation of the target spectral mean quantity $H_n(\psi)$ from (\ref{eq:tphi}) as
\begin{equation}
\label{eq:HPB-noscale}
   Q^{*}_{FDWB,n}(\psi) \equiv  n^{1/2} \{  (2\pi)^2 n^{-1}\sum_{\bj\in\J_n}\psi(\bomega_{\bj,n}) \widehat{f}_{n}(\bomega_{\bj,n})(U_{\bj}^{*}-1)\}
\end{equation}
\noindent by drawing  $ U^{*}_{\bj}$  as i.i.d exponential random variables for indices $\bj\equiv (j_1,j_2) \in \mathcal{J}_n$  involving $j_1>0$ or involving $j_1=0$ with $j_2>0$     and then setting $ U^{*}_{-\bj}\equiv U^*_{\bj}$ for the remaining indices $\bj \in \mathcal{J}_n$ (i.e., where   $j_1<0$ or where $j_1=0$ with $j_2<0$). Note that the bootstrap approximation  in (\ref{eq:HPB-noscale}) 
corresponds to the FDWB of  \cite{ng2021} for spatial data and also represents a basic step in any FDB with time series  
 (cf.~\citealp{meyer2020}). This bootstrap form aims to produce  bootstrap versions   $\{ U_{\bj}^*\widehat{f}_n(\bomega_{\bj,n})\}$ of the periodogram ordinates $\{\I_n(\bomega_{\bj,n})\}$ (i.e., indexed over $\bj\in \mathcal{J}_n$), by treating the latter as approximately independent exponential variables with respective means 
$\{f(\bomega_{\bj,n})\}$.  However, the problem in this bootstrap re-creation is that dependence among  periodogram ordinates $\I_n(\bomega_{\bj,n})$, $\bj\in\mathcal{J}_n$,  is generally not completely ignorable  so  the bootstrap quantity $Q^{*}_{FDWB,n}(\psi)$ cannot correctly capture the true variance $\sigma^2 \equiv \sigma_1^2 + \sigma_2^2$ of $H_n(\psi)$ in (\ref{eq:tphi}).  Essentially, by independently resampling or re-constructing   spatial periodogram values by  $\{ U_{\bj}^*\widehat{f}_n(\bomega_{\bj,n})\}$, the resulting FDWB in (\ref{eq:HPB-noscale}) then estimates the second variance component $\sigma_2^2$ to be zero, where the latter  may not hold in general for the underlying spatial process (i.e., unless the   process is  assumed to be Gaussian).

 To overcome the above shortcoming in FDWB, a scaling adjustment with spatial subsampling can be made to define a hybrid FDB (HFDB) rendition of $H_n(\psi)$ as
\begin{equation}
\label{eq:HPB}
H^{*}_{HFDB,n}(\psi) \equiv [\Var_{*}\{ Q^{*}_{FDWB,n}(\psi)\}+ \widehat{\sigma}^2_{2,n}]^{1/2}\frac{Q^{*}_{FDWB,n}(\psi)}{[\Var_{*}\{ Q^{*}_{FDWB,n}(\psi)\}]^{1/2}},
\end{equation}
where we first re-scale (\ref{eq:HPB-noscale}) 
to have unit variance at the bootstrap level   
using the bootstrap variance of  $\{Q^{*}_{FDWB,n}(\psi)\}$  given by 
\begin{equation*}
    \Var_{*}\{ Q^{*}_{FDWB,n}(\psi)\} = n^{-1} (4\pi^2)^2 \sum_{\bj\in \J_{n}} \psi(\bomega_{\bj,n}) (\psi(\bomega_{\bj,n}) + \psi(\bomega_{-\bj,n})) \widehat{f}_n^2(\bomega_{\bj,n})
\end{equation*} 
and then introduce a correction   
$\Var_{*}\{ Q^{*}_{FDWB,n}(\psi)\}+ \widehat{\sigma}^2_{2,n}$ to estimate the correct target variance $\sigma^2 \equiv \sigma_1^2 + \sigma_2^2$;  this scaling correction   
$\Var_{*}\{ Q^{*}_{FDWB,n}(\psi)\}+ \widehat{\sigma}^2_{2,n}$  combines 
the original bootstrap variance $\Var_{*}\{ Q^{*}_{FDWB,n}(\psi)\}$   as an approximation of the first variance component $\sigma^2_1$, with  a subsampling variance estimator
 $ \widehat{\sigma}^2_{2,n} = \widehat{\sigma}^2_n-\widehat{\sigma}^2_{1,n}$
from (\ref{eq:estvar_second}) that  approximates the second component $\sigma^2_2$. In this manner,  the hybrid bootstrap approximation $H^{*}_{HFDB,n}(\psi)$ in (\ref{eq:HPB})  can capture the correct spread in addition to the shape of the target sampling distribution of 
 $H_n(\psi)$ for spectral inference.  In the HFDB method, note that spatial subsampling plays a role in specifically estimating the variance component $\sigma^2_2$ that would otherwise be missed by FDWB, while  
 the bootstrap variance  $\Var_{*}\{ Q^{*}_{FDWB,n}(\psi)\}$ from FDWB continues to estimate the first component $\sigma_1^2$ as implicit in (\ref{eq:HPB-noscale}).  
  The next result provides a broad theoretical justification for the HFDB method.

\begin{thm}
\label{thm:hpb-consistency}
  Suppose assumptions  of Theorem \ref{thm:consistency-variance}  along with   (\ref{eq:specden-consistency}). Then, for the HFDB version $H^{*}_{HFDB,n}(\psi)$  of $H_n(\psi)$,
 \vskip0.7em
\noindent (a) HFDB spread    is consistent for the limit variance $\sigma^2 \equiv \sigma_1^2+\sigma_2^2$ of $H_n(\psi)$:
$$\Var_{*}\{H^{*}_{HFDB,n}(\psi)\} =  \Var_{*}\{Q^{*}_{FDWB,n}(\psi)\}+ \widehat{\sigma}^2_{2,n} \overset{P}{\longrightarrow} \sigma^2\ \textnormal{as}\hspace{0.2cm} n \rightarrow \infty.
$$
 \vskip0.7em
 \noindent (b) the HFDB approximation is consistent for the target distribution of $H_n(\psi)$:
\begin{eqnarray*}
     \sup_{x\in \mathbb{R}}{\left|P_{*}(H_{HFDB,n}^{*}(\psi)\leq x)-P(H_{n}(\psi)\leq x)\right|}&\overset{P}{\longrightarrow} 0& \textnormal{as} \hspace{0.2cm} n\rightarrow \infty,
\end{eqnarray*}
where $P_{*}(\cdot)$ denotes bootstrap probability   induced by   resampling.
\end{thm}

Because the HFDB  approach combines bootstrap with spatial subsampling,  
Theorem \ref{thm:hpb-consistency} shows that HFDB achieves consistency 
for distributional approximations of spectral means statistics under less  stringent moment assumptions  than previous spatial bootstraps (cf.~\citealp{ng2021}).  Furthermore, the methodology presented above for spatial HFDB can then be applied to both Gaussian and non-Gaussian processes, while  the original FDWB can only be applied to Gaussian spatial data.  Numerical   results    in Section \ref{sec:sim} further demonstrate that the HFDB can perform well, including cases where FDWB is inadequate.  

\begin{remark}
\label{rem:bias}
    As a further observation,  we can decompose the target quantity $H_n(\psi)$ as
    \[
    H_n(\psi) \equiv n^{1/2}\{\widehat{M}_{n}(\psi)-M(\psi)\}= n^{1/2}\{\widehat{M}_{n}(\psi)- E(\widehat{M}_{n}(\psi))\} + n^{1/2}\{E(\widehat{M}_{n}(\psi))-M(\psi)\},
    \] 
    where  the first part in the decomposition is a stochastic quantity with mean zero  while the second part can be treated as a non-stochastic bias.  While the above bias  decreases  to zero with increasing spatial sample size $n$,  this bias term can potentially impact approximations in small sample sizes. Therefore,  in small samples, it is also possible to consider spatial subsampling  for estimating this bias part as $\widehat{\mathrm{Bias}}_{\mathrm{sub}}\equiv L^{-1} \sum_{\ell=1}^{L} b_n^{1/2} (\widehat{M}_{\mathrm{sub}}^{(\ell)}(\psi)- \widehat{M}_n(\psi))$.  In which case, we can use $H_{HFDB,n}^{*} + \widehat{\mathrm{Bias}}_{\mathrm{sub}}$  as a bootstrap  approximation to $H_n(\psi)$, which uses subsampling to adjust bootstrap approximations  $Q^{*}_{FDWB,n}(\psi) $ from (\ref{eq:HPB-noscale}) for both   variance as in (\ref{eq:HPB}) as well as for bias.  We numerically examine such bias adjustments in  Section~\ref{sec:sim}. 
\end{remark}

\subsection{Implementation of HFDB in practice}
\label{sec:HFDB-inpractice}

Implementing HFDB requires selecting both a block/subsample size and a kernel-based spectral density estimator in a principled manner.
For block sizes, 
one may consider blocks $\mathcal{D}_b =[1,n_1^{(b)}]\times [1,n_2^{(b)}]$ defined by
scaling $n_k^{(b)} 
\approx C n^{1/4}=C(n_1 n_2)^{1/4}$, $k=1,2$, for some $C>0$ (often with $C=1,1.5$, or 2).  Such choice corresponds to an optimal size order for subsampling variance estimation (\citealp{nordman2004}).  This rule of thumb is fast to implement and is illustrated in numerical studies of Section~\ref{sec:sim}.   
Other related, but more computationally demanding,  approaches for block selection include the
 minimum volatility principle  of \cite{politis1999} and
the non-parametric plug-in (NPPI) method
from \cite{lahiri2007nonparametric}. The first approach  involves computing  bootstrap intervals over a range of block sizes and choosing a block around which the intervals do not change substantially according to some criterion (e.g., length). NPPI estimates optimal block scalings $n_1^{(b)}=n_2^{(b)}$ for variance estimation with a similar formula $\hat{C} n^{1/4}$ as above, through the constant $\hat{C}$ is estimated through 
additional steps   requiring  Jackknife-After-Bootstrap.

  For estimating the spectral density $f(\omega)$, $\bomega \in \Pi^2$, a suitable bivariate kernel function (e.g., Gaussian, Epanechnikov) may be used together with appropriate bandwidths $(\delta_1,\delta_2)$,   where \cite{ng2021} suggest choosing a bandwidth of the form $\delta_1=\delta_2$ with $\delta_k \in [0.05, 0.15]$ for small to moderate sample sizes. We   standardly use $\delta_1=\delta_2=0.05, 0.1$, or $0.15$  with a standard bivariate Gaussian kernel in Section~\ref{sec:sim}, which performed well with HFDB  and  also corresponds to the scheme for spectral density estimation
  with FDWB  (\citealp{ng2021}).

Once block sizes and kernel bandwidths are specified, HFDB can be used to generate bootstrap replicates for calibrating hypothesis tests and confidence intervals in the frequency domain. We recommend using at least 500 bootstrap replications, as done in Section~\ref{sec:sim}.

\section{Numerical Studies}
\label{sec:sim}

Here we examine the finite sample performance of the HFDB approximation of spatial spectral statistics over several inference problems, with comparisons to the FDWB approximation of \cite{ng2021}   (cf.~Section~\ref{sec:HPB}).   Theory from Sections~\ref{sec:sub-freq}-\ref{sec:HPB} suggests that HFDB should apply regardless of whether the underlying process is Gaussian or not (i.e.,   even when the variance component $\sigma^{2}_{2}$ from (\ref{eq:tphi}) is not zero as can occur in practice for non-Gaussian data). In cases   where  $\sigma^{2}_{2}$ is exactly zero (i.e., Gaussian processes),  HFDB   is anticipated to perform similarly to FDWB.  In the following, we present numerical results for scenarios including   Gaussian processes as well as a non-Gaussian processes. For each spatial sample size and process type  in the following,   we performed 1000  simulation runs to evaluate coverages and 500 replications  per simulated dataset to approximate  bootstrap distributions.  
For both FDWB/HFDB,    spectral density estimation 
 uses a Gaussian kernel    with  bandwidth $(\delta_1, \delta_2)=(0.05, 0.05)$,
 as in Section~\ref{sec:HFDB-inpractice}  and \cite{ng2021}, unless otherwise stated.  
Further simulation results, including additional coverage and bandwidth studies,  appear in \nameref{sec:supple}.
The HFDB method shall also implement a bias adjustment (Remark~\ref{rem:bias}) 
throughout; as shown in Section~\ref{sec:cov.coverage}, the effect of this adjustment is modest for coverage, but can offer improvements.

\subsection{Confidence intervals for spatial covariance}
\label{sec:cov.coverage}
We initially consider interval estimation of spatial covariance parameters $\gamma(\bh) \equiv \Cov(Z(\bzero),Z(\bh))$, $\bh\in\mathbb{Z}^2$,  as corresponding to a spectral mean $M(\psi)$ in (\ref{eq:spectral-mean}) with $\psi(\bomega)=\cos(\bh^{\top}\bomega)$.  
A main interest is comparing the 
 performances of HFDB and FDWB approximations
 for calibrating intervals, though we also consider alternative interval constructions through standard normal approximations of sample covariances combined with variance estimation steps.     
We focus on nominal 95\% confidence intervals, where bootstrap intervals have equi-tailed form
$[\widehat{M}_n(\psi) - \widehat{q}_{0.975} n^{-1/2}, \widehat{M}_n(\psi)-\widehat{q}_{0.025}n^{-1/2}]$
 based on quantiles $\widehat{q}_{0.025}$ and $\widehat{q}_{0.975}$ of bootstrap distributions approximated by  500 bootstrap draws for any   data set. We consider rectangular spatial datasets of sizes $30 \times 30\ (n=900), 50 \times 50\ (n=2500)$, and $70\times70\  (n=4900)$. 

\begin{figure}[!ht]
\begin{center}
  \includegraphics[height = 2.6in,width=1\textwidth]{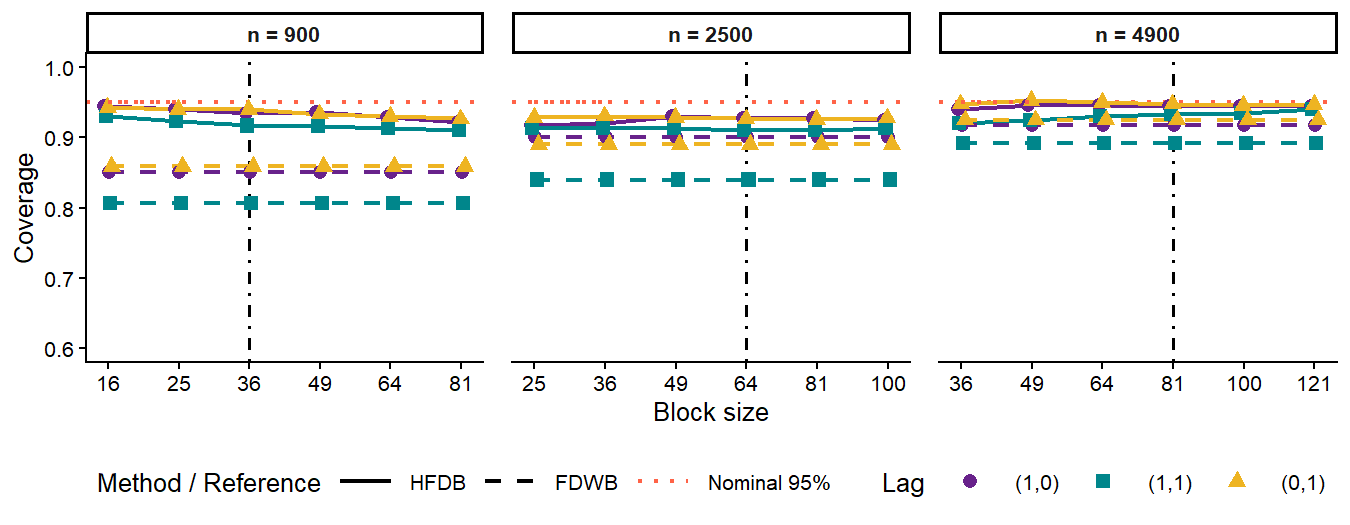}
  \caption{For covariance parameters 
  $\gamma(\bh)$ at different lags $\bh$, coverages of 95\% HFDB/FDWB intervals over various sample sizes $n$; HFDB coverages 
 are a function of  (square) block size $b_n$, where the vertical line indicates a rule of thumb block choice.}
  \label{fig1}
  \end{center}
\end{figure}

 As a first study, we consider a real-valued mean-zero   spatial Gaussian process with a covariance function is given by $ \gamma(\bh) = \theta^2 \exp \left(- (\|\bh\|/ \phi)^2 \right),$ where $\theta^2$ is the partial sill parameter and $\phi$ is the range parameter. For   simulation, we set  $\theta^2=1$ and $\phi =2$, with no nugget effect.   This  allows   evaluation of the interval methods for Gaussian data, where the latter  presents less higher-order dependence to complicate spectral inference. 
Figure~\ref{fig1} presents the coverage accuracy of 95\% two-tailed confidence intervals for the covariance parameter $\gamma(\bh)$  over multiple lags     $\bh = (1,0)^\top,  (0,1)^\top,
(1,1)^\top$; note that 
 HFDB coverages 
are plotted over various (square)  block sizes $b_n $ while FDWB coverages do not depend on blocks and so plot as straight lines.

 From Figure~\ref{fig1}, we observe that,
 as the overall spatial sample size $n$ increases,  both HFDB/FDWB methods show coverage accuracies that behave similarly and
 approach the nominal value for all lags. 
  This feature suggests that both HFDB and FDWB are valid for the Gaussian process
  (i.e., with $\sigma_{2}^{2} = 0$), as expected, where such similarity is anticipated to occur in this case because HFDB involves a scaling correction to FDWB that is   intended to be most impactful when  $\sigma^{2}_{2}$ in (\ref{eq:tphi}) is  non-zero.    Vertical lines in  Figure~\ref{fig1} also indicate that HFDB coverages by the rule  for block selection from Section~\ref{sec:HFDB-inpractice}  (with $C=1$ for illustration) appear adequate.

\newpage
As a non-Gaussian case, we generate a spatial field defined by  
$Z(\bs) \equiv \sum_{u=0}^{2} \sum_{v=0}^{2} a_{u,v} \varepsilon(s_1-u,s_2-v)
$, for $\bs\equiv (s_1,s_2)\in\mathbb{Z}^2$,  with i.i.d. innovations  $\{\varepsilon(\bs) \equiv X(\bs) - \mathbb{E}[X(\bs)]: \bs\in\mathbb{Z}^2\}$ using $X(\bs) \sim \mathrm{Gamma}(1/5,\sqrt{5})$ and
 filter coefficients 
$a_{u,v} = [A]_{u+1,v+1}/\|A\|$, $u,v\in\{0,1,2\}$ defined by entries $[A]_{u+1,v+1}$  of the matrix 
{\small \[A = \begin{pmatrix}
   1  &  0.3 &  0.1 \\
 0.5  &  0.2 &  0.1 \\
 0.25 &  0.1 & 0.05
\end{pmatrix}
\]}
and scaled by the Euclidean norm $\|A\|$ so that $\sum_{u,v} a_{u,v}^2 = 1$.   
Linear spatial processes are common in spatial statistics,
and this construction produces a strongly non-Gaussian field (e.g., skewness and excess kurtosis) with unit variance 
and a non-zero fourth-order cumulant structure. The latter feature implies that a variance component $\sigma_2^2>0$ impacts the distribution of spectral statistics, which the HFDB aims to accommodate.   In conjunction with this process,  we again consider interval estimation of   covariance parameters $\gamma(\bh)$ at several lags    $\bh = (1,0)^\top,  (0,1)^\top,
(1,1)^\top$. Here  the lag type  has a pronounced influence on the distribution  of covariance estimators targeted by bootstrap approximation, where horizontal or vertical lags  $\bh\in\{(1,0)^\top,  (0,1)^\top\}$ entail more  dependence than a diagonal lag $\bh=(1,1)^\top$,  leading to larger variance components (e.g., $\sigma_2^2$) and more difficulties in estimation.   

\begin{figure}[!ht]
\begin{center}
  \includegraphics[height = 2.6in,width=1\textwidth]{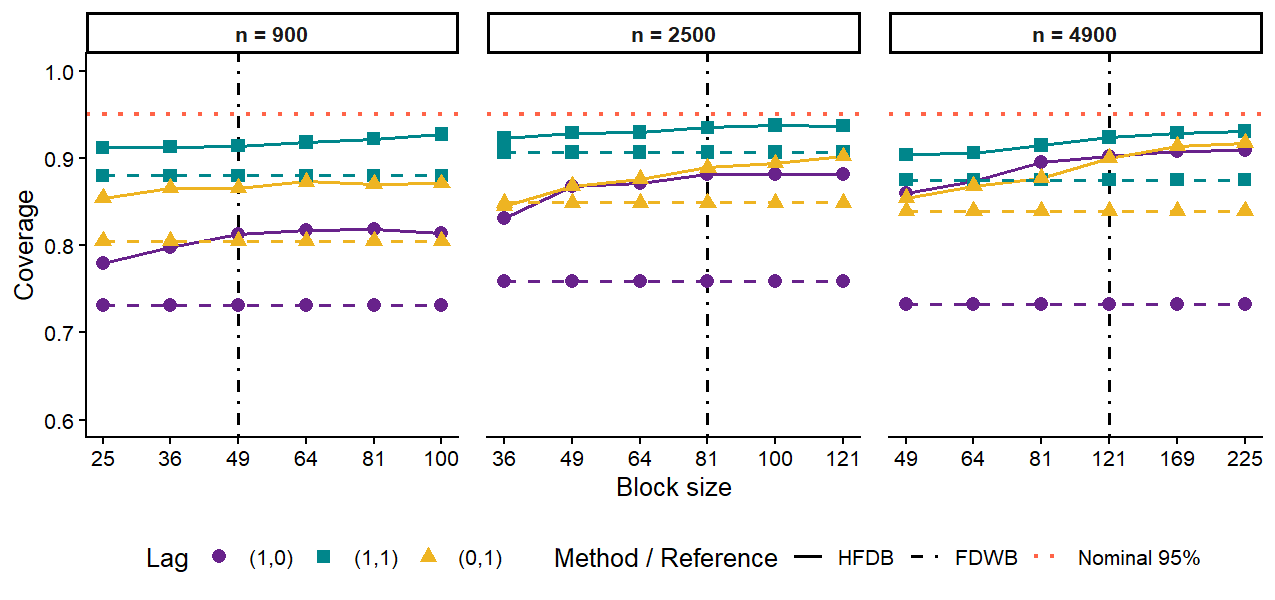}
  \caption{For covariance parameters 
  $\gamma(\bh)$ at different lags $\bh$, coverages of 95\% HFDB/FDWB intervals over various sample sizes $n$; HFDB coverages 
 are a function of  (square) block size $b_n$, where the vertical line indicates a rule of thumb block choice.}
  \label{fig2}
  \end{center}
\end{figure}

 Figure~\ref{fig2} displays the resulting coverages of HFDB/FDWB intervals, where HFDB coverages 
are again a function of   block size $b_n$.  For all covariance parameter lags,  the HFDB method consistently outperforms the FDWB in coverage across all sample sizes.
 In contrast, FDWB exhibits systematic undercoverage  due to underestimation of variance components    $\sigma_2^2$, which is most apparent with the lag $\bh=(1,0)^\top$ which induces the most dependence in estimation here (largest $\sigma_2^2$). In contrast, with $\bh = (1,1)^\top$,  the contribution of $\sigma_2^2$ is smaller, whereby underestimation of variance by FDWB  is less severe and coverage appears closer to HFDB.  The coverage performance of HFDB is relatively stable across block sizes, and 
 vertical lines in Figure~\ref{fig2} again indicate   that coverages by the block rule described in Section~\ref{sec:HFDB-inpractice} (with $C=1.25$ for illustration) seem reasonable.  

\newpage
Considering this same non-Gaussian process, we may also compare HFDB intervals to alternative (i.e., non-bootstrap) intervals based on standard normal calibrations.
In particular, we examine 95\%  intervals  $\widehat{M}_n(\psi) \pm  z_{0.975} \widehat{\sigma}_n n^{-1/2}$  
based on a normal quantile 
 $z_{0.975}=1.96$ and a subsampling variance estimator 
$\widehat{\sigma}_n$ from (\ref{eqn:sigma}) 
for the
sample covariance quantity $\widehat{M}_n(\psi)$.   

\begin{figure}[!ht]
\begin{center}
  \includegraphics[height = 2.6in,width=1\textwidth]{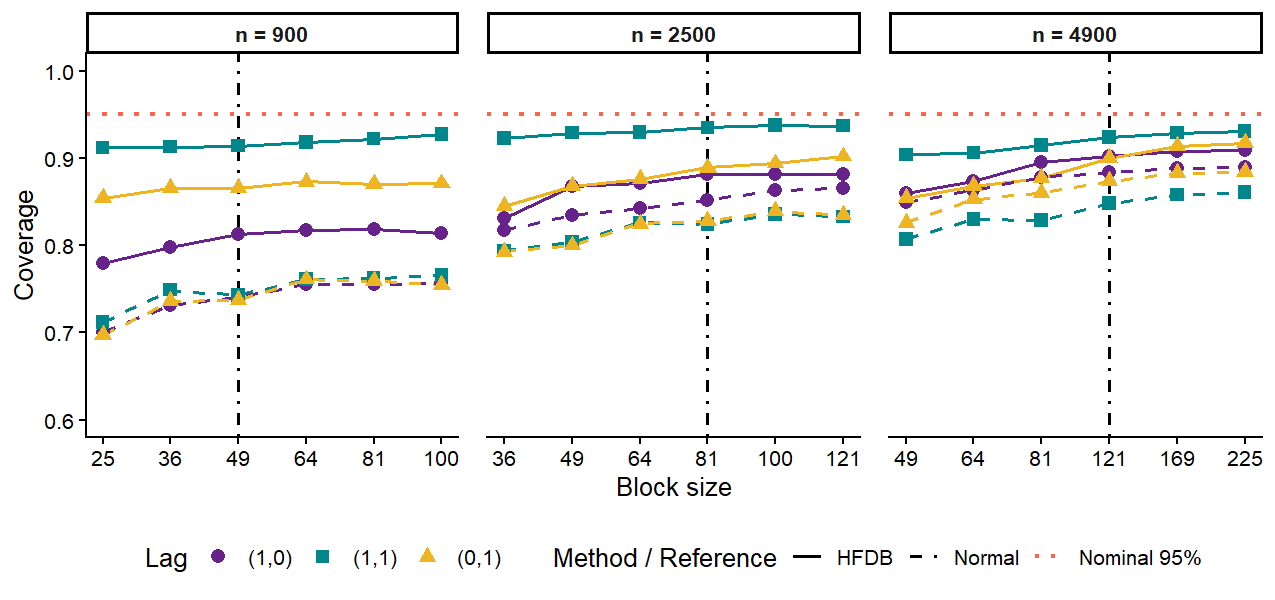}
  \caption{For covariance parameters
  $\gamma(\bh)$ at different lags $\bh$, coverages of 95\% HFDB/normal   intervals over various sample sizes $n$;  coverages of 
both intervals are a function of  (square) block size $b_n$, where the vertical line indicates a rule of thumb block choice.}
  \label{fig3}
  \end{center}
\end{figure}
Figure~\ref{fig3} provides coverage comparisons 
 between HFDB and normal-theory intervals, where  this figure serves as a direct  
  analog of  Figure~\ref{fig2} 
  after replacing FDWB intervals in the latter with normal 
  intervals.  In  Figure~\ref{fig3}, 
  the coverages of both HFDB and normal intervals depend on block size $b_n$ 
for variance estimation, though the techniques differ      between these   methods.  For all covariance parameter lags, HFDB    consistently outperforms normal intervals across block  and sample sizes; the same is not necessarily true in comparing FDWB to normal intervals across   Figures~\ref{fig2}-\ref{fig3}.  Additional simulation results in   
\nameref{sec:supple}
also show that normal intervals also  undercover 
for the previous Gaussian process as well.

As stated at the outset of Section~\ref{sec:sim}, we have focused on using a bias adjusted  version of HFDB (cf.~Remark~\ref{rem:bias}) for consistence.  The effect of this adjustment is often modest and decreasing with increasing sample sizes $n$, but can be useful.  As illustration with the non-Gaussian process,  
  Table~\ref{tab:bias-comp}
  provides a coverage  comparison of HFDB intervals
  for a covariance   
  $\gamma(\bh)$, $\bh=(1,0)^\top$, both with and without such bias adjustment.  
The adjustment provides  modest
improvements, 
  uniformly across blocks,  with less differences as  sample size  increases.

\begin{table}[!ht]
\centering
\small
\begin{tabular}{|c|cc|c|cc|c|cc|}
\hline
\multicolumn{3}{|c|}{$n = 900$} 
& \multicolumn{3}{|c|}{$n = 2500$} &\multicolumn{3}{|c|}{$n = 4900$}  \\
\hline
Block & W & W/O 
      & Block & W & W/O 
      & Block & W  & W/O \\
\hline
$5\times5$   & 78.0 & 76.2  & $6\times6$   & 83.1 & 82.5 & $7\times7$   & 86.1  & 84.9 \\
$6\times6$   & 79.8 & 76.7 & $7\times7$   & 86.8 & 84.3 & $8\times8$   & 87.4 & 86.0  \\
$7\times7\mathrlap{^{*}}$   & 81.3 & 78.7 & $8\times8$   & 87.1 & 85.3 & $9\times9$   & 89.5 & 87.2 \\
$8\times8$   & 81.7 & 78.9 & $9\times9\mathrlap{^{*}}$   & 88.2 & 85.8 & $11\times11\mathrlap{^{*}}$ & 90.2 & 88.3 \\
$9\times9$   & 81.8 & 79.5 & $10\times10$ & 88.2 & 85.9 & $13\times13$ & 90.8 & 89.0  \\
$10\times10$ & 81.4 & 79.6 & $11\times11$ & 88.1 & 86.6 & $15\times15$ & 90.9 & 89.1 \\
\hline
\end{tabular}

\caption{For $\gamma(\bh)$, $\mathbf{h} = (1,0)^\top$,  coverage (\%) of   95\%  HFDB intervals with (W) and without (W/O) bias adjustment  over   block  and sample sizes; $*$ denotes a rule of thumb block choice.}
\label{tab:bias-comp}
\end{table}

\subsection{Variogram model fitting for non-Gaussian processes}
As another illustration, we examine the coverage  performance of frequency-domain bootstrap methods for estimating a variogram model parameter. In particular, we compare FDWB with   HFDB  for the same non-Gaussian process from Section~\ref{sec:cov.coverage}.  
We consider a working semi-variogram model of the form
$
\nu(\bh;\theta) \equiv \theta   g_0(\bh)$ with $g_0(\bh) \equiv 1 - \exp\left(-\|\bh\|/\phi_0\right)$ for $\bh\in\mathbb{Z}^2$, 
where $\phi_0 =2.5$ and only the sill parameter $\theta>0$ is estimated. Given a set of lags $\{\bh_1,\dots,\bh_3\}\equiv\{(1,0), (1,1), (0,1)\}$, the parameter   $\theta$ can be formulated through the least squares criterion defined in~\ref{sec:examples}.  

  \begin{figure}[!ht]
\begin{center}
  \includegraphics[height = 2.6in,width=1\textwidth]{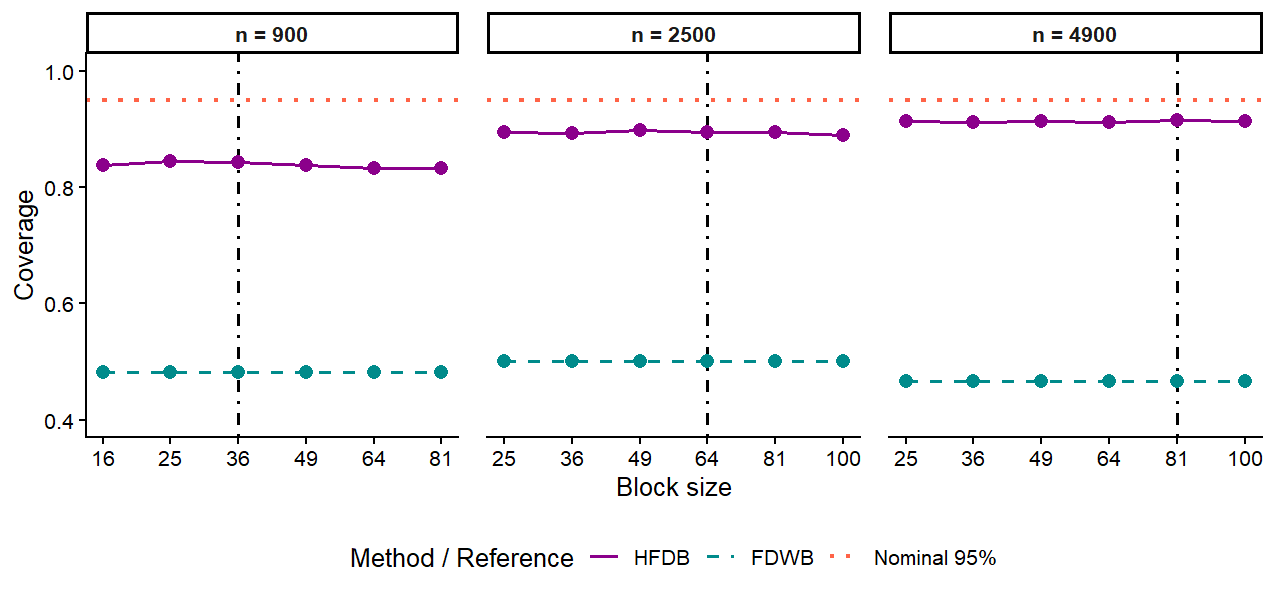}
  \caption{Coverages of 95\% HFDB/FDWB intervals for the semi-variogram parameter with different sample sizes $n$; HFDB coverages 
 are a function of  (square) block size $b_n$, where the vertical line indicates a rule of thumb block choice.}
  \label{fig4}
  \end{center}
\end{figure}

For different  sampling region sizes $n=900,2500,4900$, 
  Figure~\ref{fig4} shows the coverage rates of nominal 95\% confidence intervals for the semi-variogram parameter $\theta$ based on FDWB and HFDB methods, with the latter    coverages as a function of   block size  $b_n $.
   Here the FDWB performs poorly due to an underestimation of variance, leading to a substantial undercoverage. In contrast, HFDB significantly improves coverage which approaches the nominal level as the sample size increases, exhibiting consistent performance across block sizes.  The vertical line in  Figure~\ref{fig4} also suggests  good coverage by the block rule of thumb  ($C=1$ for illustration)  from Section~\ref{sec:HFDB-inpractice}.
    For bootstrap inference, this numerical study further demonstrates the importance of accounting for higher-order dependence that can appear in the variance of frequency domain statistics with non-Gaussian spatial processes,
    where HFDB   leads to
     more reliable inference.

\newpage
\subsection{Calibration of Spatial Isotropy Tests}
\label{sec:6}
 Here we examine the performance of bootstrap methods when applied to tests of spatial isotropy (i.e., the covariance function of the stationary spatial process is rotationally invariant). To construct a test, we  follow the setup of \cite{guan2004nonparametric}, which has also been adopted in \cite{ng2021}. In particular, we consider  testing the following hypothesis:
 $$
 H_0: 2 \kappa(\bh_i) = 2 \kappa(\bh_j),\ \mbox{for all} \hspace{0.2cm} \bh_i, \bh_j \in \Lambda, \;\bh_i \neq \bh_j,\ \mbox{and}\ ||\bh_i|| = ||\bh_j||, 
 $$
 where $\Lambda \equiv \{\bh_1, \bh_2,\ldots,\bh_m\}$ is a prespecified set of sites and $2\kappa(\bh) \equiv E(Z(\bm{0}) - Z(\bh))^2$ is the variogram at lag $\bh$.
Letting $G  \equiv (2\kappa(\bh_1),\ldots, 2\kappa(\bh_m))$ denote  the vector of variograms in $\Lambda$, 
 \cite{guan2004nonparametric} proposed a test of spatial isotropy by assessing whether
 $A G =0$ holds for a given full row rank matrix $A$ of linear contrasts.  The same test can    also be written in terms of a spectral mean assessment of whether $M(\psi) = 0$ holds for various spectral mean functionals.  
For illustration, if we consider lags $\Lambda = \{\bh_1,\bh_2\}$ for $\bh_1\equiv (1, 0)$ and $\bh_2 \equiv (0, 1)$ along $A =[-1\quad 1]$, then the test of isotropy involves  assessing whether
spectral mean statistic    $\widehat{M}_n(\psi)$ based on $\psi(\bomega) \equiv \{2 \cos(\bh_1^{\top} \bomega) - 2\cos(\bh_2^{\top} \bomega) \}$  is estimating a spectral mean $M(\psi)=0$ of zero under $H_0$.  
\cite{guan2004nonparametric} 
proposed subsampling for setting critical values for this test statistic, 
but we may also examine the effectiveness of the proposed HFDB method 
 as well as FDWB from \cite{ng2021}. 
In testing, both FDWB and  HFDB versions of this test statistic 
are computed  as in (\ref{eq:HPB-noscale})  and (\ref{eq:HPB}), respectively,  upon applying $\psi(\bomega) \equiv \{2 \cos(\bh_1^{\top} \bomega) - 2\cos(\bh_2^{\top} \bomega) \}$; these bootstrap statistics are then centered to have mean zero which mimics how the original-data test statistic $ H_{n}(\psi) \equiv \sqrt{n} \widehat{M}_n(\psi) $ behaves  under isotropy $H_0$ with $M(\psi)=0$.

For purposes of simulating spatial data to examine size control and power in testing, we use a mean zero Gaussian process with a spherical covariance function,
\begin{eqnarray}
\label{eq: spherical cov}
    \gamma(\bh) \equiv \begin{cases}
        \sigma^2 \left( 1- \frac{3r}{2\phi} + \frac{r^3}{2\phi^3}\right) + \eta \mathbb{I} \{\bh = \bm{0}\}    & \text{if} \hspace{0.2cm} 0 \leq r \leq \phi, \\
         \eta \mathbb{I} \{\bh = \bm{0}\} & \text{otherwise},
    \end{cases}
\end{eqnarray}
 where $\sigma^2$ is the partial sill parameter, $\phi$ is the range parameter, $\eta$ is the nugget
effect, and $r \equiv \sqrt{\bh^\top B \bh}$ is a distance related to a matrix $B$, described next, as geometric anisotropy transformation.
Given an anisotropy angle $\tau_A$ and anisotropy ratio $\tau_R$, define the
rotation matrix $R$ and shrinking matrix $T$ as 
\begin{eqnarray*}
    R \equiv \begin{bmatrix}
        \cos(\tau_A) & \sin(\tau_A) \\
        -\sin(\tau_A) & \cos(\tau_A)
    \end{bmatrix}
 \hspace{0.2cm} \text{and} \hspace{0.2cm} T \equiv \begin{bmatrix}
        1 & 0 \\
        0 & \tau_R
    \end{bmatrix},
\end{eqnarray*}
then we may set $B \equiv R^{'}T^{'}TR$ is a $2 \times 2$ positive definite matrix representing a geometric anisotropy transformation. A random field with spherical covariance function as in (\ref{eq: spherical cov}) is called generally anisotropic  unless $\tau_R = 1$ holds, which corresponds to the isotropic case.
We consider  model parameters $(\eta, \sigma^2, \phi, \tau_A, \tau_R) = (2, 3, 4, 0, \tau_R)$ for different anisotropy ratios $\tau_R$ to    generate   mean-zero Gaussian data $G(\cdot)$, following \cite{guan2004nonparametric} and \cite{ng2021}. 
To consider a non-Gaussian process, we generate a Gaussian random field  as above, and define observations as
$Z(\bs) = G(\bs) + 0.15 \big(G(\bs)^2 - \mathbb{E}[G(\bs)^2]\big)$.  This transformation preserves the process second-order structure, while inducing nonzero fourth-order cumulants  (i.e., $\sigma_2^2>0$) in the distribution of test statistics, which represents a situation where   differences between FDWB, HFDB, and spatial subsampling can become most apparent.

\begin{table}[!ht]
    \centering  \small
        \caption{Rejection rates of nominal 5\% spatial isotropy test (isotropy if $\tau_R=1$ here) with Subsampling (Sub), FDWB, and HFDB calibrations over  Gaussian and non-Gaussian data.}
        \label{tab:power}
        \begin{tabular}{|llll|l|lll|} 
        \multicolumn{4}{c}{Gaussian} &\multicolumn{1}{c}{}& \multicolumn{3}{c}{non-Gaussian}\\
           \cline{1-4} \cline{6-8}  
             $\tau_R$ & Sub   & FDWB & HFDB &\hspace*{1cm} &                $\tau_R$ & Sub  & HFDB \\    
            \cline{1-4}   \cline{6-8}                                                      
            $1$    & 0.030  & 0.044   & 0.044     & &              $ 1$  & 0.024  & 0.055 \\
            $1.1$  & 0.167  & 0.195   & 0.195     & &              $1.2$ & 0.149  & 0.260 \\           
            $1.2$  & 0.520  & 0.541   & 0.541     & &              $1.4$ & 0.534  & 0.683 \\           
            $1.3$  & 0.823  & 0.819    & 0.819    & &              $1.6$ & 0.857  & 0.931  \\          
            $1.4$  & 0.964  & 0.976   & 0.976     & &              $1.8$ & 0.962  & 0.989  \\          
            $1.5$  & 0.997  & 1       & 1         & &              $2$   & 0.993  & 0.997 \\          
   \cline{1-4} \cline{6-8}                                                          
        \end{tabular}                      
        \end{table}


Table~\ref{tab:power} displays empirical rejection probabilities (based on 1000 simulations) for the Gaussian and non-Gaussian processes, respectively, based on sample sizes of $50 \times 50$ with a nominal testing level of 5\%; similar results with 10\% nominal level  appear in the \nameref{sec:supple}. We have used block sizes of $7 \times 7$ for the Gaussian process, as in \cite{ng2021}, and $11 \times 11$ with bandwidth $(\delta_1, \delta_2) = (0.15,0.15)$ for illustration. Both FDWB and HFDB perform similarly for Gaussian processes, as expected. In particular, both maintain size close to the nominal level of 5\% when the null hypothesis of isotropy is true (i.e., $\tau_R=1$), and both show increasing power as $\tau_R$ deviates from 1, outperforming spatial subsampling in both size control and power. However, under non-Gaussianity,   FDWB largely fails to maintain   size   (rejection rate of $0.114$ for $\tau_R =1$)  to the extent that of being unsuitable for power comparison and is excluded in Table~\ref{tab:power} for this case.
That is, FDWB exhibits inflated rejections  and a lack of control over the Type I error. In contrast, both HFDB and spatial subsampling remain effective under non-Gaussianity, with HFDB showing some power advantages over subsampling as the anisotropy ratio increases.

\section{Concluding Remarks} 
\label{sec:conclusion}
The proposed spatial resampling method, called Hybrid Frequency Domain Bootstrap (HFDB), combines two resampling approaches in the form of subsampling and bootstrap in order to validly approximation the distribution of spatial spectral mean statistics.   
 Spectral means are fundamental spatial parameters, yet spectral statistics have complicated distributions owing to complex variances that are difficult to estimate. Previous bootstrap methods, such as the FDWB proposed by \cite{ng2021}, rely on specific distributional assumptions and can fail to provide valid distributional approximations for non-Gaussian spatial data, due to issues in correctly capturing the variances of spectral statistics for such data. In contrast, the HFDB   overcomes these limitations by incorporating a scaling adjustment from spatial subsampling to correct the spread of bootstrap approximations in the frequency domain.  This 
 leads to more accurate and general estimation of 
  the sampling distribution of spatial spectral mean statistic, as suggested by theory and numerical results in Sections~\ref{sec:sub-freq}-\ref{sec:sim}.
 
While we have focused on gridded spatial data here, the general HFDB 
approach of approximating the complicated distribution of spectral mean statistics, by    combining periodogram resampling with subsampling variance corrections,
can potentially be extended to inference about spatial processes with irregular sampling locations.  However, several new challenges arise in this setting that are not present in the gridded  
spatial case.  Firstly, due to the irregular pattern of    sampling locations, the spatial periodogram can have a multiplicative 
bias (i.e., estimating
$K f(\bomega)$ rather than the true spectral density $f(\bomega)$ for some $K\neq 1$, cf.~\citealp{bandy2010}).  As a consequence, spectral mean estimators also have multiplicative biases along with variance expressions/corrections that can change with the pattern of spatial locations.  Additionally, spectral mean estimators based on irregularly located spatial data  also require defining or specifying a grid of frequency spacings  to evaluate discrete integrals.  This aspect is less natural than for spatial lattice data and can create non-trivial additive bias for spatial spectral means  (i.e., non-vanishing bias in large samples); see \cite{subba2018} 
or \cite{zhang2024statistical}.
A modified version of HFDB could potentially accommodate such biases through bias adjustments (as initially explored here) or appropriate frequency grid choices, but requires further study.  Finally,   general distributional theory for    spectral means under non-regular locations remains a difficult and open problem.  While spectral mean statistics can have normal limits, even for stochastically located spatial data,  current results
are mainly restricted to uniform sampling locations  (cf.~\citealp{subba2018}, 
\citealp{zhang2024statistical}), whereby biases  (related to periodogram bias above) and variances have somewhat simplified expressions.    
HFDB developments for irregular spatial data could provide useful distribution approximations, though require more investigation.  

\section*{Acknowledgments}
The authors are thankful to NSF for the  NSF ACCESS MTH250064 allocation and for support by NSF DMS-2514857 (Bera and Bandyopadhyay) and 2515719 (Nordman).

\section*{Appendix}
\label{sec:appendix}
For notational simplicity, we shall use $b$ rather than $b_n$ in the following. The proofs of Theorem~\ref{thm:consistency-variance} and Theorem~\ref{thm:hpb-consistency}
involve several technical Lemmas~\ref{thm:lem1}-\ref{thm:lem3}  provided below, which are   established in  \nameref{sec:supple}.  Define $G^{*}(\bomega_{\bj,n}) \coloneqq \widehat{f}_n(\bomega_{\bj,n}) U^{*}_{\bj}$, where $ U^{*}_{\bj}$s are i.i.d standard exponential drawn using the exponential resampling mechanism.

\begin{lemma}
\label{thm:lem1}
 Suppose Assumptions~\ref{assumption:assump1}-\ref{assumption:assump3} hold with $ H_n(\psi) \equiv n^{1/2}\{\widehat{M}_{n}(\psi)-M(\psi)\} \overset{D}{\rightarrow} \mathcal{N}(0, \sigma^2) $ and  $\widehat{M}_{n}(\psi), M(\psi)$,  $\sigma^2\equiv \sigma^2(\psi)$ from  Section~\ref{sec:spectral mean parameter}-\ref{sec:spec-distn}. Then, if subsample size $b \equiv n_1^{(b)}n_2^{(b)}$ satisfies $1/n_k^{(b)} + n_k^{-1} n_k^{(b)} \rightarrow 0$ for $k=1,2$ as $n \equiv n_1n_2\rightarrow \infty$,\\
$(i)$\, $\sup_{x\in \mathbb{R}}{\left|L^{-1} \sum_{l=1}^{L} \mathbb{I}\{b^{1/2}(\widehat{M}^{(\ell)}_{\mathrm{sub}}(\psi) - M(\psi)) \le x\}-P(H_{n}(\psi)\leq x)\right|}\overset{P}{\rightarrow} 0$;\\
   $(ii)$ \, $\sup_{x\in \mathbb{R}}{\left|L^{-1} \sum_{l=1}^{L} \mathbb{I}\{b^{1/2}(\widehat{M}^{(\ell)}_{\mathrm{sub}}(\psi) - \widetilde{M}_n(\psi)) \le x\}-P(H_{n}(\psi)\leq x)\right|}\overset{P}{\rightarrow} 0$;\\
    $(iii)$ \, $b^{1/2}(\widetilde{M}_n(\psi)- M(\psi)) \overset{P}{\rightarrow} 0 $; \quad
       $(iv)$ \, $L^{-1}\sum_{l=1}^{L} b(\widehat{M}^{(\ell)}_{\mathrm{sub}}(\psi) - M(\psi))^2 \overset{P}{\rightarrow} \sigma^2(\psi) $;\\
and        $(v)$\, $L^{-1}\sum_{l=1}^{L} b(\widehat{M}^{(\ell)}_{\mathrm{sub}}(\psi) - \widetilde{M}_n(\psi))^2 \overset{P}{\rightarrow} \sigma^2(\psi) $, \\
with subsample quantities $\widehat{M}^{(\ell)}_{\mathrm{sub}}(\psi)$,  $\widetilde{M}_n(\psi)$ from Section~\ref{sec:SFDS}.
\end{lemma}

\begin{lemma}
 \label{thm:lem2}   
 Suppose Assumptions \ref{assumption:assump1}-\ref{assumption:assump4} hold
 with  
 $\psi: \Pi^2 \rightarrow \mathbf{R}$ of bounded variation. Then, if the subsample size $b\equiv n_1^{(b)}n_2^{(b)}$ satisfies $1/n_k^{(b)} + n_k^{-1} n_k^{(b)} \rightarrow 0$ for $k=1,2$ as $n\equiv n_1 n_2 \rightarrow \infty$,
 \begin{eqnarray*}
     4\pi^{2} b^{-1} \sum_{\bj \in \J_b} \psi(\bomega_{\bj,b}) L^{-1} \sum_{\ell=1}^{L} \left(\I_{\mathrm{sub}}^{(\ell)}(\bomega_{\bj,b}) - \widetilde{\I}_n(\bomega_{\bj,b})\right)^2 \overset{P}{\rightarrow} \int_{\Pi^2} \psi(\bomega) f^2(\bomega) d\bomega .
 \end{eqnarray*}
\end{lemma}

\begin{lemma}
\label{thm:lem3}   
Under Lemma~\ref{thm:lem2} conditions with 
 $Q^{*}_{FDWB,n}(\psi)$ from Section~\ref{sec:HPB},  as $n\to \infty$,\\
   $(i)$ $\Var_{*}\{ Q^{*}_{FDWB,n}(\psi)\} \overset{P}{\rightarrow} \sigma^2_1(\psi)$;
   $(ii) $ $\sup_{x\in \mathbb{R}} \left| P_{*}(Q^{*}_{FDWB,n}(\psi) \le x) - \Phi(x/\sigma_1(\psi)) \right|\overset{P}{\rightarrow} 0$,
where $\Var_{*}$, $P_{*}$   denote  variance and probability induced by resampling, and  $\Phi(\cdot)$ is the standard normal distribution function.
\end{lemma}
\subsubsection*{Proof of Theorem \ref{thm:consistency-variance}} 
Using Lemma~\ref{thm:lem2}, we have $\widehat{\sigma}^2_{1,n}(\psi) \overset{P}{\rightarrow} \sigma^2_1(\psi)$.
Since, $\psi(\cdot)$ is of bounded variation, so is
$\psi(\bomega)(\psi(\bomega)+\psi(-\bomega))$. Using Slutsky's theorem and Lemma~\ref{thm:lem1}(v) we have
$\widehat{\sigma}^2_{2,n}(\psi)= \widehat{\sigma}^2_n(\psi) - \widehat{\sigma}^2_{1,n}(\psi) \overset{P}{\rightarrow} \sigma^2(\psi) - \sigma^2_1(\psi) = \sigma^2_{2}(\psi)$. The statement of Theorem~\ref{thm:consistency-variance} now follows. \qedsymbol

\subsubsection*{Proof of Theorem \ref{thm:hpb-consistency}}
Part~(a) follows by applying Lemma~\ref{thm:lem3}(i) and Slutsky's theorem to give $\Var_{*}\{Q^{*}_{FDWB,n}(\psi)\}+ \widehat{\sigma}^2_{2,n}(\psi) \overset{P}{\rightarrow} \sigma^2(\psi)$.   To  prove the consistency of the HFDB distribution in part~(b), consider an arbitrary sequence $\{n_m\}$ of $\{n\}$. Then,  from Lemma~\ref{thm:lem3} and Theorem~\ref{thm:consistency-variance},  there exists a further subsequence $\{n_{m_k}\}$ of $\{n_m\}$ such that almost sure convergence holds as
\begin{eqnarray*}
    \sup_{x\in \mathbb{R}} \left| P_{*}(Q^{*}_{FDWB,n_{m_k}}(\psi) \le x) - \Phi(x/\sigma_1(\psi)) \right| \overset{a.s.}{\rightarrow} 0; \hskip1em
    \Var_{*}\{ Q^{*}_{FDWB,n_{m_k}}(\psi)\} \overset{a.s.}{\rightarrow} \sigma^2_1(\psi);
\end{eqnarray*}
and 
$\Var_{*}\{Q^{*}_{FDWB,n_{m_k}}(\psi)\}+ \widehat{\sigma}^2_{2,n_{m_k}}(\psi) \overset{a.s}{\rightarrow} \sigma^2(\psi)$ as $n_{m_k} \rightarrow \infty$.
Using these facts with
the HFDB definition,  
along the subsequence $\{n_{m_k}\}$, we have $H^{*}_{HFDB,n_{m_k}}(\psi) \overset{D}{\rightarrow} \mathcal{N}(0, \sigma^2(\psi))$ 
 by  Slutsky's theorem,  or equivalently  
 $\sup_{x\in \mathbb{R}} | P_{*}(H^{*}_{HFDB,n_{m_k}}(\psi) \le x) - \Phi(x/\sigma(\psi))  | \overset{a.s}{\rightarrow}0$ as $n_{m_k} \rightarrow \infty$. Since the choice of the sub-sequence was arbitrary, we have
\begin{eqnarray*}
     \sup_{x\in \mathbb{R}} \left| P_{*}(H^{*}_{HFDB,n}(\psi) \le x) - \Phi(x/\sigma(\psi)) \right| \overset{P}{\rightarrow} 0 \hspace{0.2cm}\textnormal{as} \hspace{0.2cm} n \rightarrow \infty .
\end{eqnarray*}
Now using $\sup_{x\in \mathbb{R}} \left| P(H_{n}(\psi) \le x) - \Phi(x/\sigma(\psi)) \right| = o_P(1)$  completes the proof. \qedsymbol

\bibliographystyle{chicago}
\bibliography{bibliographySFDS}

\clearpage

\phantomsection
\label{sec:supple}

\begin{center}
    {\LARGE\bfseries Supplementary Material}
\end{center}

\renewcommand{\thesection}{S\arabic{section}}
\renewcommand{\thetable}{S\arabic{table}}
\renewcommand{\thefigure}{S\arabic{figure}}
\renewcommand{\theequation}{S\arabic{equation}}

\setcounter{section}{0}
\setcounter{table}{0}
\setcounter{figure}{0}
\setcounter{equation}{0}

\begin{abstract}
These supplementary materials consist of additional numerical results as well as further mathematical proofs, which are presented in Sections~\ref{sec:s1} and \ref{sec:s2}, respectively.  
In particular, Section~\ref{sec:s1}
contains further simulation
evidence regarding the performance of spatial bootstrap methods (FDWB/HFDB) 
for spectral 
inference (Section~\ref{sec:s1.1}),  additional numerical details on tests of spatial isotropy (Section~\ref{sec:s1.2}), as well as   simulation studies regarding kernel bandwidths for bootstrap methods (Section~\ref{sec:s1.3}).  Section~\ref{sec:s2} provides proofs of technical lemmas needed to justify the proposed bootstrap procedure (HFDB) of the main manuscript.  
\end{abstract}

\section{Numerical Studies}
\label{sec:s1}

To provide additional numerical evidence in conjunction with the simulation studies presented in Section~5 of the main manuscript, here we further investigate the performance of the FDWB and HFDB methods for inferring spectral means for both Gaussian and non-Gaussian spatial processes, where non-Gaussianity is introduced through transformations. This comparison aims to highlight the coverage accuracy of the methods under different settings.

\subsection{Confidence intervals for spatial covariance}
\label{sec:s1.1}

Continuing with the example in Section~5.1, we consider a real-valued, mean-zero weakly stationary spatial process with a Gaussian covariance function. The parametric form of the covariance function is given by $ \gamma(\bh) = \theta^2 \exp \left(- (\|\bh\|/ \phi)^2 \right),$ where $\theta^2$ is the partial sill parameter and $\phi$ is the range parameter. For the simulations, we set  $\theta^2=1$ and $\phi =2$, with no nugget effect. We consider rectangular spatial datasets of sizes $30 \times 30\ (n=900), 50 \times 50\ (n=2500)$, and $70\times70\  (n=4900)$. We report results for nominal 90\% confidence intervals, where bootstrap intervals have equi-tailed form
$[\widehat{M}_n(\psi) - \widehat{q}_{0.95} n^{-1/2}, \widehat{M}_n(\psi)-\widehat{q}_{0.05}n^{-1/2}]$
 based on quantiles $\widehat{q}_{0.05}$ and $\widehat{q}_{0.95}$ of bootstrap distributions approximated by  500 bootstrap draws for each data set.

\begin{figure}[ht]
\begin{center}
  \includegraphics[height = 2.58in,width=1\textwidth]{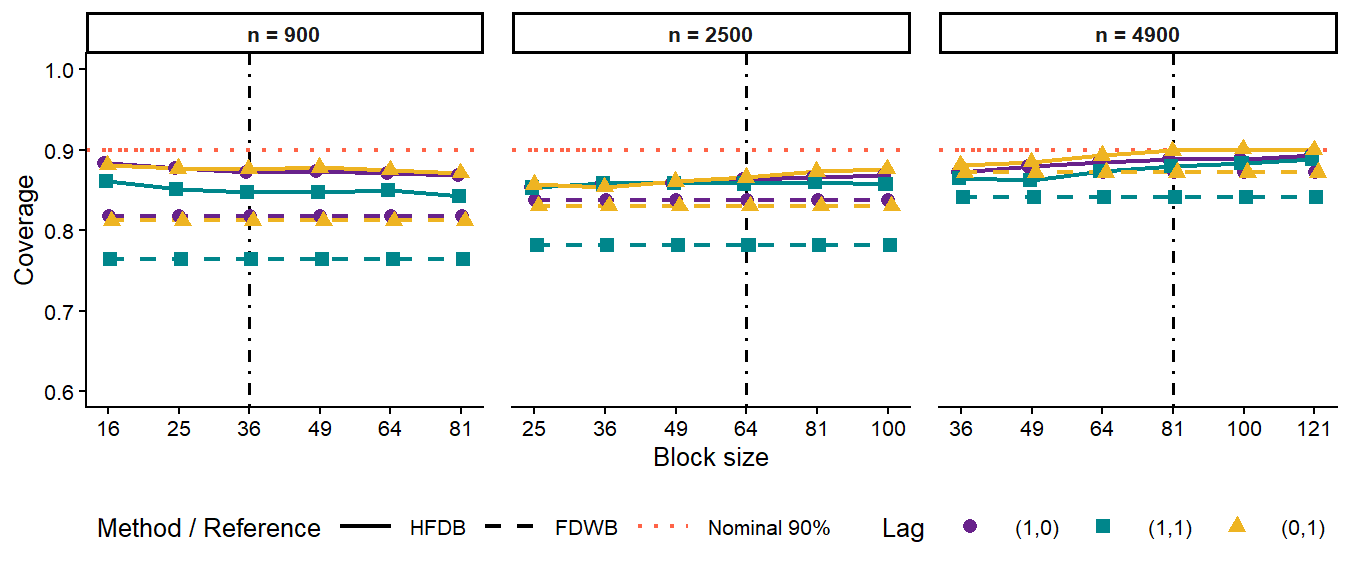}
  \caption{Empirical coverages of 90\% HFDB/FDWB intervals for the covariance parameter $\gamma(\bh)$, based on different subsample sizes $b_n$, lags $\bh$, and sample sizes $n$. The vertical line indicates a rule of thumb block choice.}
  \label{sup-fig1}
  \end{center}
\end{figure}

Figure~\ref{sup-fig1} presents the coverage accuracy of 90\% two-tailed confidence intervals for the covariance parameter. HFDB performs better compared to FDWB overall. For small sample sizes, both methods show some undercoverage, especially FDWB. As the spatial sample size $n$ increases, their coverage improves and gets closer to the nominal level. However, for $\bh=(1,1)^\top$, FDWB still shows some undercoverage even at large sample sizes, while HFDB reaches the nominal level. Overall, HFDB and FDWB behave similarly for this Gaussian process (where $\sigma_2=0$ holds in the sampling distributions targeted by bootstrap), but HFDB performs slightly better. This behavior is consistent with the fact that the scaling correction in HFDB becomes most relevant only when  $\sigma_2^2$ in $\sigma^2$ is non-zero. The coverage performance of HFDB is relatively stable across block sizes, and the vertical lines in Figure~\ref{sup-fig1} indicate the coverages corresponding to the block rule described in Section~4.1 (with $C=1$ for illustration), which appear reasonable.

Subsequently,  we present a further comparison of HFDB and the normal approximation for confidence interval construction for the same Gaussian process using the same target parameter. As discussed in Section~3.1, we can estimate the variance of spectral mean using Eq.~3.3 via subsampling variance estimator and create confidence intervals using normal approximations. For comparison, we construct 95\% confidence intervals using the normal approximation, together with a normal quantile $\text{z}_{0.975} =1.96$, yielding intervals of the form $\widehat{M}_n(\psi) \pm \text{z}_{0.975} \times \widehat{\sigma}_n n^{-1/2}$.

\begin{figure}[!ht]
\begin{center}
  \includegraphics[height = 2.6in,width=1\textwidth]{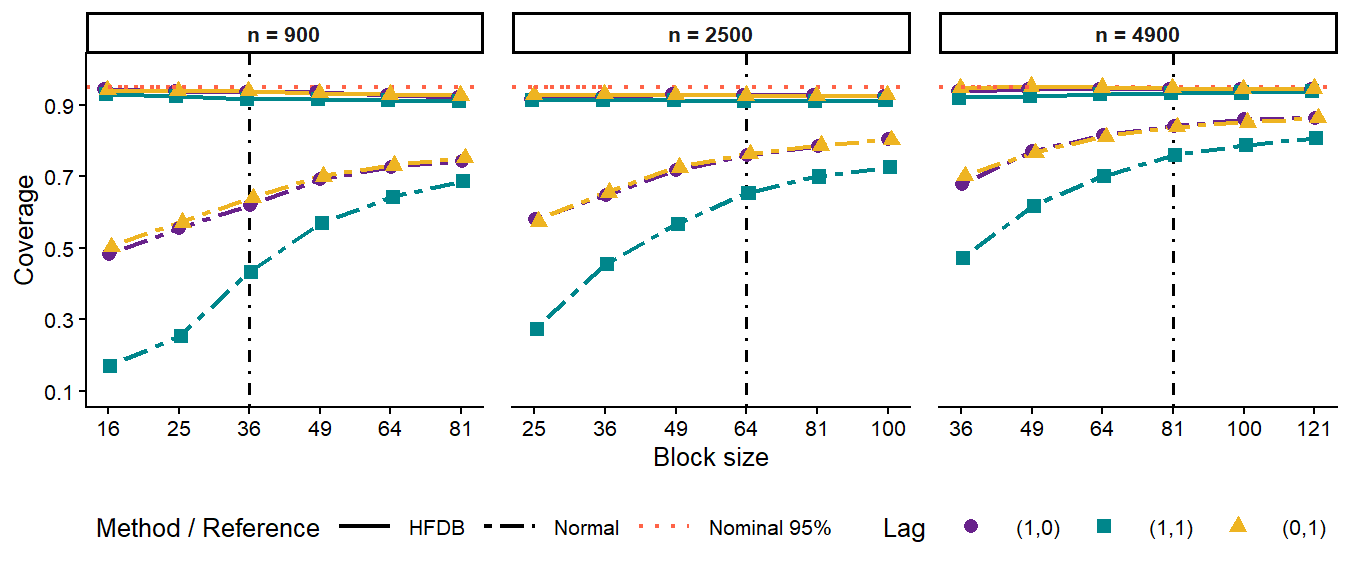}
  \caption{For covariance parameters 
  $\gamma(\bh)$ at different lags $\bh$, coverages of 95\% HFDB/normal approximation intervals over various sample sizes $n$; both HFDB/normal coverages 
 are a function of  (square) block size $b_n$, where the vertical line indicates a rule of thumb block choice.}
  \label{sup-fig2}
  \end{center}
\end{figure}

From Figure~\ref{sup-fig2}, we observe that HFDB consistently outperforms normal approximation  across block sizes and sample sizes, even in the case of a Gaussian process. HFDB maintains close to the nominal coverage for all lags $\bh$  and across all sample sizes  $n$. In contrast, the normal approximation generally undercovers, especially for smaller sample sizes,   though these normal intervals improve and gradually approach the nominal level as the block size and sample size increases. The coverage performance of HFDB is relatively stable across block sizes, and the vertical lines in Figure~\ref{sup-fig2} indicate the coverages corresponding to the block rule described in Section~4.1 (with $C=1$ for illustration).

As for another non-Gaussian example, we consider a nonlinear transformation of a stationary Gaussian spatial field. Let $\{X(\bs): \bs \in \mathbb{Z}^2\}$ be a mean-zero, variance-one stationary Gaussian process with Mat\'ern-type covariance function (cf.~\citealp{stein1999}). We generate the non-Gaussian field $\{Z(\bs): \bs \in \mathbb{Z}^2\}$ according to
\[
Z(\bs) = X(\bs) + 0.3\{X(\bs)^2 - 1\}, \quad \bs \in \mathbb{Z}^2.
\]
For the Gaussian field, we set  $\nu=1$, \(\alpha =1/7\), and \(\phi= 0.007\) to ensure that the process variance is  close to 1.
This model is reasonable, though challenging, for our study because it preserves spatial stationarity while introducing non-Gaussian behavior through a  nonlinear transformation. In particular, the process exhibits nonzero fourth-order cumulants, making the second variance component $\sigma_2^2$ relevant, when trying to provide interval estimators of covariance parameters  $
\gamma(\bh) \equiv \operatorname{Cov}\{Z(\bs), Z(\bs+\bh)\} 
$ over different lags $\bh$.
 Using the same simulation setup for a Gaussian process $X(\bs)$   as described in Section~5 of the main manuscript, 
we examine coverages of 95\% FDWB and HFDB intervals over different spatial sample sizes $30 \times 30, 50 \times 50,$ and $70 \times 70$. and for the three lags $\bh=(1,0)^\top$, $(1,1)^\top$, and $(2,0)^\top$.  

\begin{figure}[ht]
\begin{center}
  \includegraphics[height = 2.6in,width=1\textwidth]{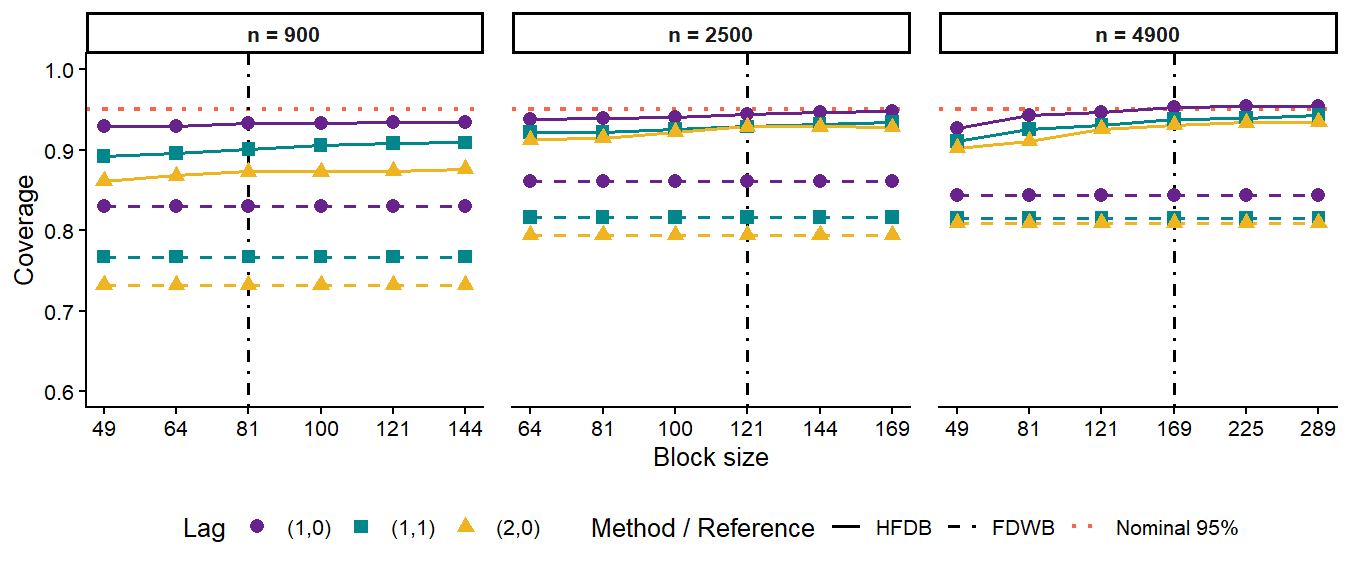}
  \caption{Empirical coverages of 95\% HFDB/FDWB intervals for the covariance parameter $\gamma(\bh)$, based on different subsample sizes $b_n$, lags $\bh$, and sample sizes $n$. The vertical line indicates a rule of thumb block choice.}
  \label{sup-fig3}
  \end{center}
\end{figure}

We observe from the coverage results shown in  Figure~\ref{sup-fig3} that FDWB underperforms across all three lags, performing worse at higher lags,  with only marginal improvement as the sample size increases. In contrast, HFDB achieves better coverage, with performance gradually improving toward the nominal level as the sample size increases for all three lags. The coverage performance of HFDB is relatively stable across
block sizes, and the vertical lines in Figure~\ref{sup-fig3} indicate the coverages corresponding to the
block rule described in Section~4.1 (with $C = 1.5$ for illustration).

\subsection{Calibration of Spatial Isotropy Tests}
\label{sec:s1.2}

As described in Section 5.3 of the main manuscript, the simulation results for a sample region $50 \times 50$ , using the same parameters as before for both Gaussian and non-Gaussian processes, are based on 1000 simulations with 500 replications at a 10\% nominal level and are as follows.

\begin{table*}[!ht]
    \centering
    \begin{minipage}[t]{0.48\textwidth}
        \centering
        \caption{Rejection rates for Gaussian data.}
        \label{sup-table1}
        \begin{tabular}{|llll|}
            \hline
            $\tau_R$ & Sub & FDWB & HFDB \\ 
            \hline
            $1$    & 0.083  & 0.091  & 0.091  \\
            $1.1$  & 0.281  & 0.295   & 0.295  \\
            $1.2$  & 0.665  & 0.656  & 0.656  \\
            $1.3$  & 0.910  & 0.903  & 0.903  \\
            $1.4$  & 0.990  & 0.991  & 0.991  \\
            $1.5$  & 0.999  & 1      & 1 \\
            \hline
        \end{tabular}
    \end{minipage}
    \hfill
    \begin{minipage}[t]{0.48\textwidth}
        \centering
        \caption{Rejection rates for Non-Gaussian data.\textsuperscript{a}}
        \label{sup-table2}
        \begin{tabular}{|lll|}
            \hline
            $\tau_R$ & Sub & HFDB \\ 
            \hline
            $ 1$  & 0.077   & 0.114 \\
            $1.2$ & 0.332   & 0.375 \\
            $1.4$ & 0.754   & 0.787 \\
            $1.6$ & 0.945   & 0.955  \\
            $1.8$ & 0.994   & 0.992  \\
            $2$   & 0.997   & 1 \\
         \hline
        \end{tabular}
    \end{minipage}
    
    \vspace{0.5em}
    \captionsetup{format=plain}
       \caption*{\footnotesize\textsuperscript{a} For non-Gaussian data, FDWB fails to maintain the size at the 10\% level, with a rejection rate of $0.174$ for $\tau_R =1$, indicating its unsuitability for such data.}
\end{table*}

Tables~\ref{sup-table1}-\ref{sup-table2} indicate that both frequency domain resampling methods perform well when the underlying distribution is Gaussian, as expected. Both HFDB and FDWB perform better in terms of size and power than spatial subsampling (Sub) when the underlying process is Gaussian. However, under non-Gaussianity, FDWB fails to maintain size under isotropy, making it unsuitable for such scenarios due to inflated rejection rates and poor control of the Type I error. In contrast, both HFDB and subsampling remain effective, with HFDB outperforming spatial subsampling (Sub) in both size and power.

\subsection{Sensitivity to Bandwidth Selection}
\label{sec:s1.3}

  Our objective here is to examine the sensitivity of the bootstrap methods (FDWB/HFDB) to bandwidth selection in kernel spectral density estimation.  Recall, from Section~4, that both FDWB/HFDB as  frequency domain bootstraps require such density estimation.  As outlined in Section~4.1, 
the HFDB approach uses a bivariate standard Gaussian kernel with a joint bandwidth $(\delta_1,\delta_2)$  such that $\delta_1=\delta_2 \in [0.05,0.15]$ as suggested by \cite{ng2021} for FDWB in moderate sample sizes.

\begin{figure}[ht]
\centering
\includegraphics[height = 2.13in,width=0.6\textwidth]{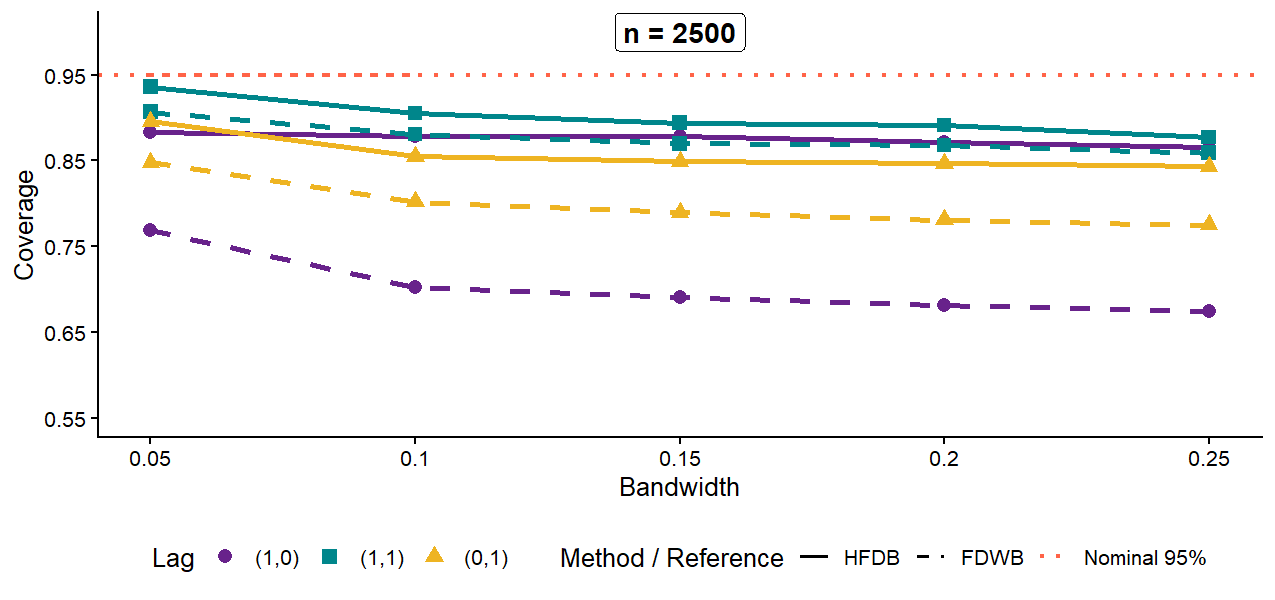}
\caption{Coverage performance of 95\% HFDB/FDWB intervals under different bandwidth choices.}
\label{fig:bw-hfdb}
\end{figure}

As a  study on the effect  of bandwidth $\delta_1=\delta_2$ choice for both bootstrap methods (FDWB/HFDB),
 Figure~\ref{fig:bw-hfdb} shows the resulting coverage of 95\%  intervals 
for covariance parameters
$\gamma(\bh)$, $\bh =(1,0)^\top, (1,1)^\top$, and $(0,1)^\top$, with spatial data simulated from the non-Gaussian process of Section~5.1 on size $n=2500$ 
sample region (i.e., $50\times 50$).  For illustration, we use a bandwidth range $\delta_1 \in [0.05,0.25]$ that is somewhat larger than the range $[0.05,0.15]$ proposed by \cite{ng2021}. Additionally, the HFDB
method uses a $10\times 10$ block size based on the rule of thumb block choice from Section~4.1 (e.g., $C=1.4$).

From Figure~\ref{fig:bw-hfdb}, both bootstraps 
(FDWB/HFDB)   exhibit coverages that similarly    decrease  with increasing bandwidths here, so that the shorter bandwidth range $[0.05,0.15]$ of \cite{ng2021} appears reasonable.   While bandwidth can effect coverage, and the best coverages in  Figure~\ref{fig:bw-hfdb} appear to occur similarly for both bootstraps at  $\delta_1=0.05$, the overall coverage of HFDB appears better across lags and also HFDB appears less sensitive than FDWB to bandwidth choice 
at some lags (e.g., $\bh =(1,0)^\top$).

\section{Proofs of Technical Lemmas}
\label{sec:s2}
Here we establish the technical lemmas, mentioned in the Appendix of the main manuscript, as needed for establishing the proposed subsampling variance estimators and bootstrap distributional estimators (i.e., Theorems~3.1 and 4.1 of the main paper).

\subsection{Proof of Lemma~1}
To formalize our proofs, we introduce the following: for $r=1,2$, and $x \in \mathbb{R}$, \\
$H_n \equiv H_n(\psi)$;\hspace{0.3cm} $Y \sim \mathcal{N}(0, \sigma^2(\psi))$,\hspace{0.3cm} $t_n^{(r)} \equiv E(H_n^r)$, \hspace{0.3cm} $t_{n,C}^{(r)} \equiv E\left( H_n -E(H_n)\right)^r$, \\
$F_n(x) \equiv P(H_n \le x)$, \hskip1em$F(x) \equiv P(Y \le x)$,\\
$\widehat{F}_n(x) \equiv L^{-1} \sum_{\ell=1}^{L} \mathbb{I} [ b_n^{1/2} (\widehat{M}^{(\ell)}_{\mathrm{sub}}(\psi) - M(\psi)) \le x ]$; 
$ \widehat{F}_{n,C}(x) \equiv L^{-1}  \sum_{\ell=1}^{L} \mathbb{I} [ b_n^{1/2} ( \widehat{M}^{(\ell)}_{\mathrm{sub}}(\psi) - \widetilde{M}_n(\psi)) \le x ]$;\\
$\widehat{t}_n^{(r)} \equiv L^{-1} \sum_{\ell=1}^{L}  \left\{b_n^{1/2} \left( \widehat{M}^{(\ell)}_{\mathrm{sub}}(\psi) - M(\psi)\right) \right\}^{r}$; $ \widehat{t}_{n,C}^{(r)} \equiv L^{-1}\sum_{\ell=1}^{L} \left\{ b_n^{1/2} \left( \widehat{M}^{(\ell)}_{\mathrm{sub}}(\psi) - \widetilde{M}_n(\psi)\right) \right\}^r$. \\

These definitions establish key quantities used to prove the asymptotic validity of bootstrap procedures. They distinguish between population-level targets (e.g., $M(\psi)$) and sample-based estimates (e.g., $\widetilde{M}_n(\psi)$) to carefully track bias and variability in the resampling framework. We define $r$-th order central moments and empirical distributions centered around both the true and estimated targets.

\noindent First we prove statement 1(i). By notation defined above, we have to show $\sup_{x\in \mathbb{R}} |\widehat{F}_n(x)- F(x)|\overset{P}{\rightarrow} 0$ as $n \rightarrow \infty$. Let $\mathcal{C} \subset \mathbb{R}$ is the set of continuity points of $F$. It suffices to show $\widehat{F}_n(x_0) \overset{P}{\rightarrow} F(x_0)$ as $n \rightarrow \infty$ for an arbitrary fixed $x_0 \in \mathcal{C}$. From this, one may take a countable, dense(in $\mathbb{R}$) collection $\{x_i: i \ge 1 \} \subset \mathcal{C}$ and, for any sub-sequence  $\{n_m\}$ of $\{n\}$, extract a further sub-sequence $\{n_{m_k}\} \equiv \{n_k\}$ where, almost surely(a.s.), it holds that $\widehat{F}_{n_{m_k}}(x_i) \rightarrow F(x_i)$ as $n_{m_k} \rightarrow \infty$ for each $i \ge 1$; the latter implies that $\sup_{x\in \mathbb{R}}|\widehat{F}_{n_{m_k}}(x)- F(x)|\overset{a.s.}{\rightarrow} 0$, which is equivalent to $\sup_{x\in \mathbb{R}} |\widehat{F}_n(x)- F(x)|\overset{P}{\rightarrow} 0$ as $\{n_m\}$ was arbitrary. Then statement (i) follows from a triangle inequality. 

\vskip1em

\noindent Now, due to boundedness,  $\widehat{F}_n(x_0)\overset{P}{\rightarrow} F(x_0)$ is equivalent to $E(\widehat{F}_n(x_0)- F(x_0))^2 \rightarrow 0$, and we show the latter.  

\noindent To control this, fix $\lambda \in (0,1)$ where
$\lambda$ sets a threshold for spatial proximity, and for $\bj=(j_1,j_2), \bk=(k_1,k_2)$ integer vectors in the plane (the $i$-th component can be $0,1,...,\left(n_i-n^{(b)}_i\right)$ for $i=1,2$) we split the double sum of covariances in the MSE into indices $\bj, \bk$ such that:

\begin{itemize}
  \item[(a)]  $|j_i - k_i| < \lambda n_i$ for some $i = 1,2$,
  \item[(b)] and those for which $|j_i - k_i| \geq \lambda n_i$ for both $i = 1,2$.
  
\end{itemize}

In case (a), the number of such dependent index pairs has an upper bound $n\lambda / L$.

In case (b), the assumption
\[
\frac{|j_i - k_i|}{n^{(b)}_i} \to \infty \quad \text{for} \quad i = 1,2,
\]
along with $n^{(b)}_i / n_i \to 0$ (or $n_i / n^{(b)}_i \to \infty$ ), implies that the covariance between $\mathbb{I}(H_{\bj,b} \leq x_0)$ and $\mathbb{I}(H_{\bk,b} \leq x_0)$ must disappear asymptotically. That is,
\begin{eqnarray*}
&&\sup_{\substack{|j_i - k_i| \geq \lambda n_i \\ \forall i=1,2}} 
\text{Cov}\left(\mathbb{I}(H_{\bj,b} \leq x_0), \mathbb{I}(H_{\bk,b} \leq x_0)\right) \\
&&= \sup_{\substack{|j_i - k_i| \geq \lambda n_i \\ \forall i=1,2}}
\left|P(H_{\bj,b} \leq x_0, H_{\bk,b} \leq x_0) - 
P(H_{\bj,b} \leq x_0) P(H_{\bk,b} \leq x_0)\right|
\to 0,
\end{eqnarray*}

by Assumption 3 This guarantees
\[
D_{1,n}(\lambda) := \sup_{\substack{|j_i - k_i| \geq \lambda n_i \\ \forall i=1,2}} |P(H_{\bj,b} \leq x_0, H_{\bk,b} \leq x_0) - F(x_0)^2| \to 0,
\]
\[
D_{2,n}(\lambda) := \sup_{\substack{j_i \in \{0,1,..,(n_i-n^{(b)}_i)\} \\ \forall i=1,2}} |P(H_{\bj,b} \leq x_0) \left(F(x_0) - P(H_{\bj,b} \leq x_0)\right)| \to 0.
\]

Combining, we obtain:
\[
E \left( \widehat{F}_n(x_0) - F(x_0) \right)^2 \leq \frac{n\lambda}{L} + D_{1,n}(\lambda) + D_{2,n}(\lambda) \to \lambda.
\]

Since $\lambda > 0$ is arbitrary, we conclude:
\[
E \left( \widehat{F}_n(x_0) - F(x_0) \right)^2 \to 0,
\]
completing the proof of Lemma 1(i).

\noindent Next, we prove statements (iii) - (v) together. First we would like to prove the following statements: if $\lim_{n \rightarrow \infty } E(H_{\bj^{(n)},n}^2) = E (Y^2)$ hold for any vector of integer sequence $\bj^{(n)}$, then
\begin{eqnarray}
\label{eq:t-convergence}
    \widehat{t}_n^{(1)} \overset{P}{\rightarrow} t_n^{(1)},\hspace{0.2cm} \widehat{t}_n^{(2)} \overset{P}{\rightarrow} t_n^{(2)},\hspace{0.2cm} \textnormal{and}  \hspace{0.2cm} \widehat{t}_{n,C}^{(2)} \overset{P}{\rightarrow} t_{n,C}^{(2)} \hspace{0.2cm}\textnormal{as} \hspace{0.2cm} n \rightarrow \infty.
\end{eqnarray}
We show $\widehat{t}_n^{(r)} \overset{P}{\rightarrow} t_n^{(r)}$, for $r=1,2$. We would also have $\widehat{t}_{n,C}^{(2)} \overset{P}{\rightarrow} t_{n,C}^{(2)}$ by Slutsky's theorem. \\
 For any $C > 0$ such that $\ensuremath{\pm}C \in \mathcal{C}$, it follows from statement $(i)$ that
\begin{equation}
\label{eq:tCon-inP}
  \widehat{t}_n^{(r)}(C) \equiv \frac{1}{L} \sum_{\ell =1}^{L} H_{\bj^{(\ell)},b}^{r} \mathbb{I} \{ |H_{\bj^{(\ell)},b}| \le C\} = \int y^r \mathbb{I} \{ |y| \le C\} d\widehat{F}_n(y) 
  \overset{P}{\rightarrow} E(Y^r \mathbb{I}\{|Y| \le C\}).
 \end{equation}
That is, by statement (i), for any sub-sequence $\{n_j\}$ of $\{n\}$, there exists again further sub-sequence $\{n_k\} \subset \{n_j\}$ where $\sup_{x\in \mathbb{R}} |\widehat{F}_{n_k}(x)- F(x)|\overset{a.s.}{\rightarrow} 0$ as $\{n_k\} \rightarrow 0$. Then using Continuous Mapping theorem, for a random variable $Y_{n_k}^{*} \sim \widehat{F}_{n_k}$, we have $(Y_{n_k}^{*})^r \mathbb{I}\{|Y_{n_k}^{*}| \le C\} \overset{D}{\rightarrow} Y^r \mathbb{I}\{ |Y \le C\}$, and so the corresponding expected values (bounded) converge as 
$\int y^r \mathbb{I} \{ |y| \le C\} d\widehat{F}_{n_k}(y) 
  \overset{a.s.}{\rightarrow} E(Y^r \mathbb{I}\{|Y| \le C\})$ as $n_k \rightarrow \infty$ implying Eq.~(\ref{eq:tCon-inP}). For given $\upsilon > 0$ and $\lambda \in (0,1)$, pick and fix $C$ such that $ E \left(Y^r \mathbb{I}\{|Y |> C\} \right) < \upsilon\lambda/3$. Then it suffices to show 
  \begin{equation}
  \label{eq:limsup}
     \limsup_{n \rightarrow \infty} P\left(|\widehat{t}_n^{(r)}- \widehat{t}_{n}^{(r)}(C)| > \upsilon/3 \right) \le \lambda,
  \end{equation}
from which
\begin{eqnarray*} 
\limsup_{n \rightarrow \infty} P\left(|\widehat{t}_n^{(r)}- E(Y^r)| > \upsilon \right) 
&\le& \limsup_{n \rightarrow \infty} P\left( |\widehat{t}_n^{(r)}(C)- E(Y^r \mathbb{I}\{|Y | \le C\})| > \upsilon/3 \right)\\
&&\hskip1em +\limsup_{n \rightarrow \infty} P\left(|\widehat{t}_n^{(r)}- \widehat{t}_{n}^{(r)}(C)| > \upsilon/3 \right) \le \lambda    
\end{eqnarray*}
follows by Eq.~(\ref{eq:tCon-inP})-(\ref{eq:limsup}) and $P\left(E(Y^r \mathbb{I}\{|Y| > C\}) > \upsilon/3 \right) = 0$, which gives us $\widehat{t}_n^{(r)} \overset{P}{\rightarrow} E(Y^r)$. Since
we have assumed for any vector of integer sequence $\bj^{(n)}$, $E(H_{\bj^{(n)},n}^2) \rightarrow E(Y^2)$ as $n \rightarrow \infty$, we have
$t_n^{(2)} \equiv E(H_n^2) \rightarrow E(Y^2)= \sigma^2(\psi)$ and, by dominated convergence theorem (DCT), $E(|H_n|) \rightarrow E(|H|)$ and $t_n^{(1)} \equiv E(H_n) \rightarrow E(Y)=0$ as $n \rightarrow \infty$. Thus, we have established $\widehat{t}_n^{(r)} \overset{P}{\rightarrow} t_n^{(r)}$ for $r=1,2$. 

\vskip1em
\noindent By Markov's inequality, for $r=1,2$, we have 
$$
    P\left(|\widehat{t}_n^{(r)}- \widehat{t}_n^{(r)}(C)| > \upsilon/3 \right) \le \frac{\upsilon}{3} \max_{1 \le \ell \le L} E\left(|H_{\bj^{(\ell)},b}|^r \mathbb{I} \{|H_{\bj^{(\ell)},b}|> C\}\right),
$$
so that, Eq.~(\ref{eq:limsup}) follows from $\limsup_{n \rightarrow \infty} \max_{1 \le \ell \le L} E \left( |H_{\bj^{(\ell)}, b}|^r \mathbb{I}\{|H_{\bj^{(\ell)},b}| > C\}\right) \le \upsilon\lambda/3$ which is implied by the moment conditions. Note that a contrary result 
$$
\limsup_{n \rightarrow \infty} \max_{1 \le \ell \le L} E \left( |H_{\bj^{(\ell)}, b}|^r \mathbb{I}\{|H_{\bj^{(\ell)},b}| > C\}\right) > \upsilon\lambda/3
$$ would imply existence of a sub-sequence $\{n_j\}$ and a positive integer $m_{n_j}$ such that 
$$
E \left( |H_{\bj^{(m_{n_j})},b_{n_j}}|^r \mathbb{I}\{|H_{\bj^{(m_{n_j})},b_{n_j}}| > C\} \right) > \upsilon \lambda/3
$$ holds for each $n_j$, which contradicts
$$
    \lim_{n_j \rightarrow \infty} E \left( |H_{\bj^{(m_{n_j})},b_{n_j}}|^r \mathbb{I}\{|H_{\bj^{(m_{n_j})},b_{n_j}}| > C\} \right) = E \left(|Y|^r \mathbb{I}\{|Y| > C\}\right) < \upsilon\lambda/3,
$$
that follows by DCT from $H_{\bj^{(m_{n_j})}, b_{n_j}} \overset{D}{\rightarrow} Y $ along with $E(|H_{\bj^{(m_{n_j})}, b_{n_j}}|^r ) \rightarrow E(|Y|^r)$ from the assumed moment conditions. 

\vskip1em
\noindent Now, the remaining is to show $\lim_{n \rightarrow \infty}E\left( H^2_{\bj^{(n)},n}\right) = E(Y^2)$. By 4th-order stationarity we only need to prove $\lim_{ \rightarrow \infty} E(H_n^2)= \sigma^2(\psi)$. By Assumption 1, \citealp{brillinger2001}, and \citealp{fuentes2002} we have,
\begin{eqnarray}
\label{eq:periodogram mean}
E(\I_n(\bomega))= f(\bomega) + O(n^{-1}), \hspace{0.3cm} \bomega \in \Pi^2  
\end{eqnarray}
and for $\bomega_1, \bomega_2 \in \Pi^2$, and $i=1,2$,
\begin{eqnarray}
\label{eq:cov periodogram}
\textnormal{Cov}(\I_n(\bomega_1), \I_n(\bomega_2)) = \begin{cases} 
  f^2(\bomega_1) + O( n^{-1})  & \text{if } |\omega_{1i}|=|\omega_{2i}| (\neq 0, \hspace{0.2em} \forall i=1,2) \\
  \frac{(2\pi)^2}{n}f_4(\bomega_1, \bomega_2, -\bomega_2) + O( n^{-1}) & \text{if } |\omega_{1i}|\neq|\omega_{2i}|.
\end{cases}
\end{eqnarray}
The error terms are uniform in $\bomega$. Hence, for $i=1,2$,
\begin{eqnarray*}
\textnormal{Var}(H_n(\psi)) &=& n^{-1}(4\pi^2)^2\sum_{\bj, \bk\in \J_n} \psi(\bomega_{\bj, n}) \psi(\bomega_{\bk, n}) \textnormal{Cov} (\I_n(\bomega_{\bj,n}), \I_n(\bomega_{\bk,n})) \\
&=& n^{-1}(4\pi^2)^2\left[\sum_{\bj \in \J_n} \psi(\bomega_{\bj, n})\left( \psi(\bomega_{\bj, n}) +  \psi(-\bomega_{\bj, n})\right) \textnormal{Var} (\I_n(\bomega_{\bj,n})) \right.\\
&&\hskip1em \left.+ \sum\limits_{\substack{\bj,\bk \in \J_n \\ |\bomega_{j_i}| \neq |\bomega_{k_i}|}} \psi(\bomega_{\bj, n}) \psi(\bomega_{\bk, n}) \textnormal{Cov} (\I_n(\bomega_{\bj,n}), \I_n(\bomega_{\bk,n})) \right]\\
&=:& S_1+S_2.
\end{eqnarray*}

By Eq.~ (\ref{eq:cov periodogram}) and Riemann sum form, it holds that
\begin{eqnarray*}
    S_1 &=& (2\pi)^2 n^{-1} 4\pi^2 \sum_{\bj \in \J_n}\psi(\bomega_{\bj, n})\left( \psi(\bomega_{\bj, n}) +  \psi(-\bomega_{\bj, n})\right) f^2(\bomega_{\bj,n}) + 
    O( n^{-1}) \\
    &&= (2\pi)^2 \int_{\Pi^2} \psi(\bomega)(\psi(\bomega)+\psi(-\bomega)) f^2(\bomega) d\bomega + O( n^{-1}) = \sigma^2_1(\psi) + O( n^{-1}).
\end{eqnarray*}
Similarly, we have that
\begin{eqnarray*}
    S_2 &=& n^{-1}(4\pi^2)^2 \sum\limits_{\substack{\bj,\bk \in \J_n \\ |\bomega_{j_i}| \neq |\bomega_{k_i}|}} \psi(\bomega_{\bj, n}) \psi(\bomega_{\bk, n}) \left( \frac{(2\pi)^2}{n}f_4(\bomega_{\bj, n}, \bomega_{\bk, n}, -\bomega_{\bk, n}) + O( n^{-1})\right) \\
    &=& (2\pi)^2 n^{-2} 4\pi^2 \left[\sum_{\bj,\bk \in \J_n} \psi(\bomega_{\bj, n}) \psi(\bomega_{\bk, n})f_4(\bomega_{\bj, n}, \bomega_{\bk, n}, -\bomega_{\bk, n}) \right.\\
    &&\hskip1em \left.- \sum_{\bj \in \J_n} \psi(\bomega_{\bj, n})(\psi(\bomega_{\bj, n})+\psi(-\bomega_{\bj, n}))f_4(\bomega_{\bj, n}, \bomega_{\bj, n}, -\bomega_{\bj, n}) + O(n^{-1})\right] \\
    &=& (2\pi)^2 \int_{\Pi^{2}}\int_{\Pi^{2}} \psi(\bomega_{1})\psi(\bomega_{2})f_{4}(\bomega_{1},\bomega_{2},-\bomega_{2})\,d\bomega_{1}\,d\bomega_{2} + O( n^{-1}) = \sigma^2_2(\psi) + O(n^{-1}).
\end{eqnarray*}
Now, we have $\textnormal{Var}(H_n(\psi))= \sigma_1^2(\psi)+\sigma_2^2(\psi) + O(n^{-1}) = \sigma^2(\psi) + O(n^{-1})$. By Eq.~(\ref{eq:periodogram mean}) and bounded variation of $\psi(\cdot)$, we also have 
\begin{equation*}
E(H_n(\psi))= n^{1/2}\left(\frac{(2\pi)^2 }{n} \sum_{\bj \in \J_n} \psi(\bomega_{\bj,n})\{f(\bomega_{\bj,n}) + O(n^{-1}) \}  - \int_{\Pi^2} \psi(\bomega)f(\omega) d\omega \right) = O(n^{-1/2}).
\end{equation*}
Hence, we conclude that $E(H_n(\psi))^2= \textnormal{Var}(H_n(\psi))+ E^2(H_n(\psi))= \sigma^2(\psi)+ O(n^{-1})$. 

\vskip1em
\noindent Finally, we prove statement (ii), i.e., we will show $\sup_{x \in \mathbb{R}}\left| \widehat{F}_{n,C}(x) - F_n(x)\right| \overset{P}{\rightarrow} 0$ as $n \rightarrow \infty$.
Let $Y_n^{*}$ denotes a random variable with distribution function $\widehat{F}_{n}$, then $(Y_n^{*} - \widehat{t}^{(1)}_n)$ is a random variable
with distribution function $\widehat{F}_{n,C}$. Because $\sup_{x \in \mathbb{R}} \left| \widehat{F}_n(x) - F(x) \right| \overset{P}{\rightarrow} 0$ and $\widehat{t}_n^{(1)} \overset{P}{\rightarrow} t_n^{(1)}$, and
$t_n^{(1)} \rightarrow E(Y) = 0$ guaranteed by the moment assumption, we have, for any sub-sequence $\{n_j\} \subset \{n\}$, there
exists a further sub-sequence $\{n_k\} \subset \{n_j\}$ such that $Y^{*}_{n_k} \overset{D}{\rightarrow} Y$ and $\widehat{t}^{(1)}_{n_k} \rightarrow 0$ hold as $n_k \rightarrow \infty$ almost surely
(a.s.). By Slutsky’s theorem, $(Y^{*}_{n_k} -\widehat{t}^{(1)}_{n_k}) \overset{D}{\rightarrow} Y$ follows a.s. which, because $\{n_j\}$ was arbitrary, implies
$\sup_{x \in \mathbb{R}}\left| \widehat{F}_{n,C}(x) - F(x)\right| \overset{P}{\rightarrow} 0$. Statement (ii) follows from a triangle inequality. \qedsymbol
\subsection{Proof of Lemma~2}
Note that,
\begin{eqnarray*}
&&(2\pi)^2 b^{-1} \sum_{\bj \in \J_b} \psi(\bomega_{\bj,b}) \frac{1}{L} \sum_{\ell=1}^{L} \left(\I_{\mathrm{sub}}^{(\ell)}(\bomega_{\bj,b}) - \widetilde{\I}_n(\bomega_{\bj,b})\right)^2\\ 
&& =(2\pi)^2 b^{-1} \left[\sum_{\ell=1}^{L} \psi(\bomega_{\bj,b}) \left(\frac{1}{L} \sum_{\ell=1}^{L} (\I_{\mathrm{sub}}^{(\ell)}(\bomega_{\bj,b})) ^2 - 2f^2(\bomega_{\bj,b})\right)\right. \\
&&\hskip1em\left.+\sum_{\ell=1}^{L} \psi(\bomega_{\bj,b})\left(f^2(\bomega_{\bj,b}) -(\widetilde{\I}_n(\bomega_{\bj,b}))^2 \right) + \sum_{\ell=1}^{L} \psi(\bomega_{\bj,b})f^2(\bomega_{\bj,b})\right] \\
       &&=: S_{1n}+S_{2n}+S_{3n}. 
\end{eqnarray*}
It is immediate from the Riemann sum expression that, $S_{3n} \rightarrow \int_{\Pi^2} \psi(\bomega)f^2(\bomega)$ as $n \rightarrow \infty$.
Now, proceeding in the same way as \citealp{YuThesis23}, using Cauchy-Schwarz inequality and boundedness of $\psi(\cdot)$ we have $|S_{2n}| = o(1)$. Also, using the moment conditions we have $S_{1n}\overset{P}{\rightarrow} 0$, which gives us the statement of Lemma~2. \qedsymbol

\subsection{Proof of Lemma~3}
By the exponential resampling mechanism,
we have
 $\textnormal{Cov}_{*}(U_{\bj}^{*}, U_{\bk}^{*}) = \mathbb{I}\{ \bj=-\bk\}$. Hence, with $\widehat{f}_n$ being an even function it follows,
 $ \textnormal{Cov}_{*}(G^{*}(\bomega_{\bj,n}), G^{*}(\bomega_{\bk,n}))
     = \widehat{f}^{2}_n(\bomega_{\bj,n}) \mathbb{I}\{\bj=-\bk\}$ and,
 \begin{eqnarray*}
\textnormal{Var}_{*}(Q^{*}_{FDWB,n}(\psi)) &=& \frac{(4\pi^2)^2}{n} \sum_{\bj,\bk \in \J_n}  \psi(\bomega_{\bj,n}) \psi(\bomega_{\bk,n}) \textnormal{Cov}_{*}(G^{*}(\bomega_{\bj,n}), G^{*}(\bomega_{\bk,n})) \\
     &=& (2\pi)^2 \sum_{\bj \in \J_n} \frac{(2\pi)^2}{n} \psi(\bomega_{\bj,n})\{\psi(\bomega_{\bj,n}) + \psi(\bomega_{-\bj,n} )\} f^2(\bomega_{\bj,n}) \\ 
     &&\hskip1em \quad + \frac{(4\pi^2)^2}{n} \sum_{\bj \in \J_n} \psi(\bomega_{\bj,n})\{\psi(\bomega_{\bj,n}) + \psi(\bomega_{-\bj,n} )\} (\widehat{f}^2_n(\bomega_{\bj,n})- f^2(\bomega_{\bj,n})) \\
     &=:& S_1 + S_2
 \end{eqnarray*}    
Since, $S_1$ is a Riemann sum, it immediately follows that $S_1 = \sigma_1^2(\psi) + O(1)$. Using $|\J_n|= O(n)$ and Eq.~(4.7) we have $S_2=o_P(1)$ which completes the proof of Lemma~3(i). 

\vskip1em

\noindent 
For notational ease, we define $\J_n^{+} \equiv \{ \bj\equiv (j_1,j_2) \in\mathcal{J}_n:  \mbox{either $j_1>0$ or $j_1=0$ and $j_2>0$}\}$.  

\noindent 
Since $\{U_{\bj}^{*}, \bj \in \J_n^{+}\}$ are i.i.d. standard exponentially distributed and $U_{\bj}^{*}\equiv U_{-\bj}^{*}$ for $\bj \in \J_n$ where $-\bj \in \J_n^{+}$, then $Q^{*}_{FDWB,n}(\psi)$ can be written as
\begin{eqnarray}
\label{eq:qhpb-positive}
  Q^{*}_{FDWB,n} =  n^{-1/2}\sum_{\bj \in \J_n^{+}} V^{*}_{\bj,n},
\end{eqnarray}
where,
$$
   V^{*}_{\bj,n} = (2\pi)^2 (\psi(\bomega_{\bj,n}) \widehat{f}_{n}(\bomega_{\bj,n})+\psi(-\bomega_{\bj,n}) \widehat{f}_{n}(-\bomega_{\bj,n}))(U_{\bj}^{*}-1), \hspace{0.2cm} \bj \in \J_n^{+}.
$$
Since, $|\J_n^{+}| \rightarrow \infty$, as $n \rightarrow \infty$, a conditional version of Lyapunov’s CLT can
be applied to Eq.~(\ref{eq:qhpb-positive}) (using the
established convergence in probability of $\textnormal{Var}_{*}(Q^{*}_{FDWB,n}(\psi))$ to a positive constant from (i), i.e., we have to find a $\upsilon > 0$ such that 
\begin{eqnarray}
\label{eq:Lyapunov condition}
    n^{-(1+\upsilon/2)} \sum_{\bj \in \J_n^{+}} E^{*} \left(| V^{*}_{\bj, n}|^{2+ \upsilon} \right) = o_P(1).
\end{eqnarray}
Choosing $\upsilon=1$, we have 
$$
    \sum_{\bj \in \J_n^{+}} |\psi(\bomega_{\bj,n}) \widehat{f}_{n}(\bomega_{\bj,n})+\psi(-\bomega_{\bj,n}) \widehat{f}_{n}(-\bomega_{\bj,n})| \le n \cdot \sup_{\bomega\in \Pi} |\psi(\bomega)| \cdot \sup_{\bomega\in \Pi} |\widehat{f}_n(\bomega)| = O_P(n),  
$$
since both $\psi$ and $f$ (due to absolutely summable autocovariance function) are bounded, and due to Eq.~(4.7). Hence, it holds that
\begin{eqnarray}
\label{eq:O_p(n)}
     \sum_{\bj \in \J_n^{+}} \left|\psi(\bomega_{\bj,n}) \widehat{f}_{n}(\bomega_{\bj,n})+\psi(-\bomega_{\bj,n}) \widehat{f}_{n}(-\bomega_{\bj,n})\right|^3 = O_P(n).
\end{eqnarray}
Since, $U_{\bj}^{*}$ are standard exponential we have $E^{*}(|U_{\bj}^{*} - 1|^3) = 12e^{-1} -2= K(\textnormal{Say})$ gives us
\begin{eqnarray*}
    \sum_{\bj \in \J_n^{+}} \mathbf{E}^{*} \left(| V^{*}_{\bj, n}|^{2+ \upsilon} \right) = (4\pi^2)^3 K  \sum_{\bj \in \J_n^{+}} \left|\psi(\bomega_{\bj,n}) \widehat{f}_{n}(\bomega_{\bj,n})+\psi(-\bomega_{\bj,n}) \widehat{f}_{n}(-\bomega_{\bj,n})\right|^3 = O_P(n),
\end{eqnarray*}
due to Eq.~(\ref{eq:O_p(n)}). This shows Lyapunov condition holds for Eq.~(\ref{eq:Lyapunov condition}) when $\upsilon=1$. Together with Lemma~3(i), the assertion in Lemma~3(ii)  follows. \qedsymbol

\end{document}